\documentclass[graybox]{svmult}

\usepackage{mathptmx}       
\usepackage{helvet}         
\usepackage{courier}        
\usepackage{type1cm}        
\usepackage{makeidx}         
\usepackage{graphicx}        
\usepackage{multicol}        
\usepackage[bottom]{footmisc}

\usepackage{amsmath, stmaryrd, latexsym} \usepackage{euscript}
\usepackage{mathrsfs}
\usepackage{amssymb}

\begin{document}

\newcounter{bnomer} \newcounter{snomer}
\newcounter{bsnomer}
\setcounter{bnomer}{0}
\renewcommand{\thesnomer}{\thebnomer.\arabic{snomer}}
\renewcommand{\thebsnomer}{\thebnomer.\arabic{bsnomer}}

\setcounter{MaxMatrixCols}{14}

\newcommand{\sect}[1]{%
\setcounter{snomer}{0}\setcounter{bsnomer}{0}
\refstepcounter{bnomer}
\par\bigskip\begin{center}\large{\textbf{\arabic{bnomer}. {#1}}}\end{center}}
\newcommand{\fakesect}{%
\setcounter{snomer}{0}\setcounter{bsnomer}{0}
\refstepcounter{bnomer}}
\newcommand{\sst}[1]{%
\refstepcounter{bsnomer}
\par\bigskip\textbf{\arabic{bnomer}.\arabic{bsnomer}. {#1}.}}
\newcommand{\fakesst}{%
\refstepcounter{bsnomer}}
\newcommand{\defi}[1]{%
\refstepcounter{snomer}
\par\medskip\textbf{Definition \arabic{bnomer}.\arabic{snomer}. }{#1}\par\medskip}
\newcommand{\theo}[2]{%
\refstepcounter{snomer}
\par\textbf{Теорема \arabic{bnomer}.\arabic{snomer}. }{#2} {\emph{#1}}\hspace{\fill}$\square$\par}
\newcommand{\mtheop}[2]{%
\refstepcounter{snomer}
\par\textbf{Theorem \arabic{bnomer}.\arabic{snomer}. }{\emph{#1}}
\par\textsc{Proof}. {#2}\hspace{\fill}$\square$\par}
\newcommand{\mcorop}[2]{%
\refstepcounter{snomer}
\par\textbf{Corollary \arabic{bnomer}.\arabic{snomer}. }{\emph{#1}}
\par\textsc{Proof}. {#2}\hspace{\fill}$\square$\par}
\newcommand{\mtheo}[1]{%
\refstepcounter{snomer}
\par\medskip\textbf{Theorem \arabic{bnomer}.\arabic{snomer}. }{\emph{#1}}\par\medskip}
\newcommand{\mlemm}[1]{%
\refstepcounter{snomer}
\par\medskip\textbf{Lemma \arabic{bnomer}.\arabic{snomer}. }{\emph{#1}}\par\medskip}
\newcommand{\mprop}[1]{%
\refstepcounter{snomer}
\par\medskip\textbf{Proposition \arabic{bnomer}.\arabic{snomer}. }{\emph{#1}}\par\medskip}
\newcommand{\theobp}[2]{%
\refstepcounter{snomer}
\par\textbf{Теорема \arabic{bnomer}.\arabic{snomer}. }{#2} {\emph{#1}}\par}
\newcommand{\theop}[2]{%
\refstepcounter{snomer}
\par\textbf{Theorem \arabic{bnomer}.\arabic{snomer}. }{\emph{#1}}
\par\textsc{Proof}. {#2}\hspace{\fill}$\square$\par}
\newcommand{\theosp}[2]{%
\refstepcounter{snomer}
\par\textbf{Теорема \arabic{bnomer}.\arabic{snomer}. }{\emph{#1}}
\par\textbf{Схема доказательства}. {#2}\hspace{\fill}$\square$\par}
\newcommand{\exam}[1]{%
\refstepcounter{snomer}
\par\medskip\textbf{Example \arabic{bnomer}.\arabic{snomer}. }{#1}\par\medskip}
\newcommand{\deno}[1]{%
\refstepcounter{snomer}
\par\textbf{Definition \arabic{bnomer}.\arabic{snomer}. }{#1}\par}
\newcommand{\post}[1]{%
\refstepcounter{snomer}
\par\textbf{Предложение \arabic{bnomer}.\arabic{snomer}. }{\emph{#1}}\hspace{\fill}$\square$\par}
\newcommand{\postp}[2]{%
\refstepcounter{snomer}
\par\medskip\textbf{Proposition \arabic{bnomer}.\arabic{snomer}. }{\emph{#1}}%
\ifhmode\par\fi\textsc{Proof}. {#2}\hspace{\fill}$\square$\par\medskip}
\newcommand{\lemm}[1]{%
\refstepcounter{snomer}
\par\textbf{Lemma \arabic{bnomer}.\arabic{snomer}. }{\emph{#1}}\hspace{\fill}$\square$\par}
\newcommand{\lemmp}[2]{%
\refstepcounter{snomer}
\par\medskip\textbf{Lemma \arabic{bnomer}.\arabic{snomer}. }{\emph{#1}}
\par\textsc{Proof}. {#2}\hspace{\fill}$\square$\par\medskip}
\newcommand{\coro}[1]{%
\refstepcounter{snomer}
\par\textbf{Corollary \arabic{bnomer}.\arabic{snomer}. }{\emph{#1}}\hspace{\fill}$\square$\par}
\newcommand{\mcoro}[1]{%
\refstepcounter{snomer}
\par\textbf{Corollary \arabic{bnomer}.\arabic{snomer}. }{\emph{#1}}\par\medskip}
\newcommand{\corop}[2]{%
\refstepcounter{snomer}
\par\textbf{Corollary \arabic{bnomer}.\arabic{snomer}. }{\emph{#1}}
\par\textsc{Proof}. {#2}\hspace{\fill}$\square$\par}
\newcommand{\nota}[1]{%
\refstepcounter{snomer}
\par\medskip\textbf{Remark \arabic{bnomer}.\arabic{snomer}. }{#1}\par\medskip}
\newcommand{\propp}[2]{%
\refstepcounter{snomer}
\par\medskip\textbf{Proposition \arabic{bnomer}.\arabic{snomer}. }{\emph{#1}}
\par\textsc{Proof}. {#2}\hspace{\fill}$\square$\par\medskip}
\newcommand{\hypo}[1]{%
\refstepcounter{snomer}
\par\medskip\textbf{Conjecture \arabic{bnomer}.\arabic{snomer}. }{\emph{#1}}\par\medskip}
\newcommand{\prop}[1]{%
\refstepcounter{snomer}
\par\textbf{Proposition \arabic{bnomer}.\arabic{snomer}. }{\emph{#1}}\hspace{\fill}$\square$\par}

\newcommand\restr[2]{{
  \left.\kern-\nulldelimiterspace 
  #1 
  \right|_{#2} 
}}

\newcommand{\Ind}[3]{%
\mathrm{Ind}_{#1}^{#2}{#3}}
\newcommand{\Res}[3]{%
\mathrm{Res}_{#1}^{#2}{#3}}
\newcommand{\epsi}{\epsilon}
\newcommand{\tri}{\triangleleft}
\newcommand{\Supp}[1]{%
\mathrm{Supp}(#1)}

\newcommand{\reg}{\mathrm{reg}}
\newcommand{\sreg}{\mathrm{sreg}}
\newcommand{\codim}{\mathrm{codim}\,}
\newcommand{\chara}{\mathrm{char}\,}
\newcommand{\rk}{\mathrm{rk}\,}
\newcommand{\chr}{\mathrm{ch}\,}
\newcommand{\Ann}{\mathrm{Ann}\,}
\newcommand{\id}{\mathrm{id}}
\newcommand{\Ad}{\mathrm{Ad}}
\newcommand{\col}{\mathrm{col}}
\newcommand{\row}{\mathrm{row}}
\newcommand{\low}{\mathrm{low}}
\newcommand{\pho}{\hphantom{\quad}\vphantom{\mid}}
\newcommand{\fho}[1]{\vphantom{\mid}\setbox0\hbox{00}\hbox to \wd0{\hss\ensuremath{#1}\hss}}
\newcommand{\wt}{\widetilde}
\newcommand{\wh}{\widehat}
\newcommand{\ad}[1]{\mathrm{ad}_{#1}}
\newcommand{\Grr}[2]{\mathrm{Gr}(#1,#2)}
\newcommand{\Gra}[1]{\mathrm{Gr}(#1)}
\newcommand{\tr}{\mathrm{tr}\,}
\newcommand{\Cent}{\mathrm{Cent}\,}
\newcommand{\GL}{\mathrm{GL}}
\newcommand{\Or}{\mathrm{O}}
\newcommand{\SO}{\mathrm{SO}}
\newcommand{\SU}{\mathrm{SU}}
\newcommand{\SL}{\mathrm{SL}}
\newcommand{\Sp}{\mathrm{Sp}}
\newcommand{\Mat}{\mathrm{Mat}}
\newcommand{\Pf}{\mathrm{Pf}}
\newcommand{\Prim}{\mathrm{Prim}\,}
\newcommand{\Ker}{\mathrm{Ker}\,}
\newcommand{\Stab}{\mathrm{Stab}}
\newcommand{\ilm}{\varinjlim}
\newcommand{\soc}[1]{\mathrm{soc}^{#1}\,}
\newcommand{\gee}{\geqslant}
\newcommand{\lee}{\leqslant}

\newcommand{\vfi}{\varphi}
\newcommand{\teta}{\vartheta}
\newcommand{\Bfi}{\Phi}
\newcommand{\Fp}{\mathbb{F}}
\newcommand{\Np}{\mathbb{N}}
\newcommand{\Rp}{\mathbb{R}}
\newcommand{\Zp}{\mathbb{Z}}
\newcommand{\Cp}{\mathbb{C}}
\newcommand{\Hp}{\mathbb{H}}
\newcommand{\ut}{\mathfrak{u}}
\newcommand{\at}{\mathfrak{a}}
\newcommand{\nt}{\mathfrak{n}}
\newcommand{\mt}{\mathfrak{m}}
\newcommand{\htt}{\mathfrak{h}}
\newcommand{\spt}{\mathfrak{sp}}
\newcommand{\sot}{\mathfrak{so}}
\newcommand{\glt}{\mathfrak{gl}}
\newcommand{\rt}{\mathfrak{r}}
\newcommand{\zt}{\mathfrak{z}}
\newcommand{\rad}{\mathfrak{rad}}
\newcommand{\bt}{\mathfrak{b}}
\newcommand{\gt}{\mathfrak{g}}
\newcommand{\vt}{\mathfrak{v}}
\newcommand{\pt}{\mathfrak{p}}
\newcommand{\Xt}{\mathfrak{X}}
\newcommand{\Po}{\mathcal{P}}
\newcommand{\Uo}{\EuScript{U}}
\newcommand{\Fo}{\EuScript{F}}
\newcommand{\Ho}{\EuScript{H}}
\newcommand{\Do}{\mathcal{D}}
\newcommand{\Eo}{\EuScript{E}}
\newcommand{\Vo}{\EuScript{V}}
\newcommand{\Iu}{\mathcal{I}}
\newcommand{\Mo}{\mathcal{M}}
\newcommand{\Nu}{\mathcal{N}}
\newcommand{\Ro}{\mathcal{R}}
\newcommand{\Co}{\mathcal{C}}
\newcommand{\Go}{\mathcal{G}}
\newcommand{\Lo}{\mathcal{L}}
\newcommand{\Ou}{\mathcal{O}}
\newcommand{\Uu}{\mathcal{U}}
\newcommand{\Au}{\mathcal{A}}
\newcommand{\Vu}{\mathcal{V}}
\newcommand{\Bu}{\mathcal{B}}
\newcommand{\Sy}{\mathcal{Z}}
\newcommand{\Sb}{\mathcal{F}}
\newcommand{\Gr}{\mathcal{G}}
\newcommand{\Fl}{\EuScript{F}\ell}
\newcommand{\rtc}[1]{C_{#1}^{\mathrm{red}}}

\title*{Real group orbits\\ on flag ind-varieties of $\SL(\infty,\Cp)$}
\author{Mikhail V. Ignatyev\and Ivan Penkov\and Joseph A. Wolf}
\institute{Mikhail V. Ignatyev \at Department of Mathematics and Mechanics, Samara State Aerospace University, Ak. Pavlova 1, 443011 Samara, Russia, \email{\texttt{mihail.ignatev@gmail.com}}
\and Ivan Penkov \at Jacobs University Bremen, Campus Ring 1, 28759 Bremen, Germany, \email{\texttt{i.penkov@jacobs-university.de}}
\and Joseph A. Wolf \at Department of Mathematics, University of California Berkeley, Berkeley California 94720-3840, USA, \email{\texttt{jawolf@math.berkeley.edu}}}

\maketitle

\abstract{We consider the complex ind-group $G=\SL(\infty,\Cp)$ and its real forms\break $G^0=\SU(\infty,\infty)$, $\SU(p,\infty)$, $\SL(\infty,\Rp)$, $\SL(\infty,\Hp)$. Our main object of study are the $G^0$-orbits on an ind-variety $G/P$ for an arbitrary splitting parabolic ind-subgroup $P\subset G$, under the assumption that the subgroups $G^0\subset G$ and $P\subset G$ are aligned in a natural way. We prove that the intersection of any $G^0$-orbit on $G/P$ with a finite-dimensional flag variety $G_n/P_n$ from a given exhaustion of $G/P$ via $G_n/P_n$ for $n\to\infty$, is a single $(G^0\cap G_n)$-orbit. We also characterize all ind-varieties $G/P$ on which there are finitely many $G^0$-orbits, and provide criteria for the existence of open and closed $G^0$-orbits on $G/P$ in the case of infinitely many $G^0$-orbits.}

\begin{center}
\begin{tabular}{p{15cm}}
\end{tabular}
\end{center}

\section{Introduction}
\fakesect

This study has its roots in linear algebra. Witt's Theorem claims that, given any two subspaces $V_1$, $V_2$ of a finite-dimensional vector space $V$ endowed with a nondegenerate bilinear or Hermitian form, the spaces $V_1$ and $V_2$ are isometric within $V$ (i.e., one is obtained from the other via an isometry of $V$) if and only if $V_1$ and $V_2$ are isometric. When $V$ is a Hermitian space, this is a statement about the orbits of the unitary group $U(V)$ on the complex grassmannian $\Grr{k}{V}$, where $k=\dim V_1=\dim V_2$. More precisely, the orbits of $U(V)$ on $\Grr{k}{V}$ are parameterized by the possible signatures of a, possibly degenerate, Hermitian form on a $k$-dimensional space of $V$.

A general theory of orbits of a real form $G^0$ of a semisimple complex Lie group $G$ on a flag variety $G/P$ was developed by the third author in \cite{Wolf1}
and \cite{Wolf3}. This theory has become a standard tool in semisimple representation
theory and complex algebraic geometry.  For automorphic forms and automorphic
cohomology we mention \cite{WellsWolf1}, \cite{FelsHuckleberryWolf1} and
\cite{Wolf5}.  For double fibration transforms and similar applications to
representation theory see \cite{WolfZierau1} and
\cite{HuckleberryWolf3}.  For the structure of real group orbits and cycle spaces with
other applications to complex algebraic geometry see, for example,
\cite{Wolf3}, \cite{Wolf4}, \cite{Wolf5}, \cite{Wolf6},
\cite{WolfZierau1}, \cite{BarletKoziarz1},
\cite{BarciniLeslieZierau1}, \cite{FelsHucklebery1}
\cite{KroetzStanton1}, \cite{KroetzStanton2}, \cite{OrstedWolf1} and
\cite{FelsHuckleberryWolf1}.  Finally,
applications to geometric quantization are indicated by
\cite{RawnsleySchmidWolf1} and \cite{SchmidWolf2}.

The purpose of the present paper is to initiate a systematic study of real group orbits on flag ind-varieties or, more precisely, on ind-varieties of generalized flags. The study of the classical simple ind-groups like $\SL(\infty,\Cp)$ arose from studying stabilization phenomena for classical algebraic groups. By now, the classical ind-groups, their Lie algebras, and their representations have grown to a separate subfield in the vast field of infinite-dimensional Lie groups and Lie algebras. In particular, it was seen in \cite{DimitrovPenkov1} that the ind-varieties $G/P$ for classical ind-groups $G$ consist of generalized flags (rather than simply of flags) which are, in general, infinite chains of subspaces subject to two delicate conditions, see Subsection~\ref{sst:gen_flags} below.

Here we restrict ourselves to the ind-group $G=\SL(\infty,\Cp)$ and its real forms $G^0$. We study $G^0$-orbits on an arbitrary ind-variety of generalized flags $G/P$, and establish several foundational results in this direction. Our setting assumes a certain alignment between the subgroups $G^0\subset G$ and $P\subset G$.

Our first result is the fact that any $G^0$-orbit in $G/P$, when intersected with a finite-dimensional flag variety $G_n/P_n$ from a given exhaustion of $G/P$ via $G_n/P_n$ for $n\to\infty$, yields a single $G_n^0$-orbit. This means that the mapping
\begin{equation*}
\text{$\{G_n^0$-orbits on $G_n/P_n\}$ $\to$ $\{G_{n+1}^0$-orbits on $G_{n+1}/P_{n+1}\}$}
\end{equation*}
is injective. Using this feature, we are able to answer the following questions.
\begin{enumerate}
\item When are there finitely many $G^0$-orbits on $G/P$?
\item When is a given $G^0$-orbit on $G/P$ closed?
\item When is a given $G^0$-orbit on $G/P$ open?
\end{enumerate}
The answers depend on the type of real form and not only on the parabolic subgroup $P\subset G$. For instance, if $P=B$ is an upper-triangular Borel ind-subgroup of $\SL(\infty,\Cp)$ ($B$ depends on a choice of an ordered basis in the natural representation of $\SL(\infty,\Cp)$), then $G/B$ has no closed $\SU(\infty,\infty)$-orbit and has no open $\SL(\infty,\Rp)$-orbit.

We see the results of this paper only as a first step in the direction of understanding the structure of $G/P$ as a $G^0$-ind-variety for all real forms of all classical ind-groups $G$ (and all splitting parabolic subgroups $P\subset G$). Substantial work lies ahead.

\section{Background}
\fakesect

In this section we review some basic facts about finite-dimensional real group orbits.  We then discuss the relevant class of infinite-dimensional Lie groups
and the corresponding real forms and flag ind-varieties.

\subsection{Finite-dimensional case}\fakesst \label{sst:finite_dim_case}Let $V$ be a finite-dimensional complex vector space. Recall that a~\emph{real structure} on $V$ is an antilinear involution $\tau$ on $V$. The set $V^0=\{v\in V\mid\tau(v)=v\}$ is a \emph{real form} of $V$, i.e., $V^0$ is a real vector subspace of $V$ such that $\dim_{\Rp}V^0=\dim_{\Cp}V$ and the $\Cp$-linear span $\langle V^0\rangle_{\Cp}$ coincides with $V$. A real form $V^0$ of $V$ defines a unique real structure $\tau$ on $V$ such that $V^0$ is the set of fixed point of~$\tau$. A~\emph{real form} of a complex finite-dimensional Lie algebra $\gt$ is a real Lie subalgebra $\gt^0$ of $\gt$ such that $\gt^0$ is a~real form of $\gt$ as a complex vector space.

Let $G$ be a complex semisimple connected algebraic group, and $G^0$ be a \emph{real form} of $G$, i.e., $G^0$ is a real closed algebraic subgroup of $G$ such that its Lie algebra $\gt^0$ is a real form of the Lie algebra $\gt$ of $G$. Let $P$ be a parabolic subgroup of $G$, and $X=G/P$ be the corresponding flag variety. The group $G^0$ naturally acts on $X$. In \cite{Wolf1} the third author proved the following facts about the $G^0$-orbit structure of $X$, see \cite[Theorems 2.6, 3.3, 3.6, Corollary 3.4]{Wolf1} (here we use the usual differentiable
manifold topology on $X$).

\mtheo{\label{theo:finite_dim_case}\\
\textup{i)} Each $G_0$-orbit is a real submanifold of $X$.\\
\textup{ii)} The number of $G^0$-orbits on $X$ is finite.\\  \textup{iii)} The union of the open $G^0$-orbits is dense in $X$.\\ \textup{iv)} There is a unique closed orbit $\Omega$ on $X$.\\ \textup{v)} The inequality $\dim_{\Rp}\Omega\geq\dim_{\Cp}X$ holds.}

Here is how this theorem relates to Witt's Theorem in the case of a Hermitian form. Let $V$ be an $n$-dimensional complex vector space and $G=\SL(V)$. Fix a nondegenerate Hermitian form $\omega$ of signature $(p,n-p)$ on the vector space $V$ and denote by $G^0=\SU(V,\omega)$ the group of all linear operators on $V$ of determinant~$1$ which preserve $\omega$. Then $G^0$ is a real form of $G$. Given $k\leq n$, the group~$G$ naturally acts on the grassmannian $X=\Grr{k}{V}$ of all $k$-dimensional complex subspaces of $V$. To each $U\in X$ one can assign its signature $(a,b,c)$, where the restricted form $\restr{\omega}{U}$ has rank $a+b$ with $a$ positive squares and $b$ negative ones, $c$ equals the dimension of the intersection of $U$ and its orthogonal complement, and $a+b+c=k$. By Witt's Theorem, two subspaces $U_1$, $U_2\in X$ belong to the same $G^0$-orbit if and only if their signatures coincide. Set $l=\min\{p,~n-p\}$. Then one can verify the following formula for the number $|X/G^0|$ of $G^0$-orbits on $X$:
\begin{equation*}
|X/G^0|=\begin{cases}
(n-k+1)(n-k+2)/2,&\text{if }n-l\leq k,\\
(l+1)(l+2)/2,&\text{if }l\leq k\leq n-l,\\
(k+1)(k+2)/2,&\text{if }k\leq l.\\
\end{cases}
\end{equation*}
Furthermore, a $G^0$-orbit of a subspace $U\in X$ is open if and only if the restriction of $\omega$ to $U$ is nondegenerate, i.e., if $c=0$. Therefore, the number of open orbits equals $\min\{k+1,l+1\}$. There is a unique closed $G^0$-orbit $\Omega$ on $X$, and it consists of all $k$-dimensional subspaces of $V$ such that $c=\min\{k,l\}$ (the condition $c=\min\{k,l\}$ maximizes the nullity of the form $\restr{\omega}{U}$ for $k$-dimensional subspaces $U\subset V$). In particular, if $k=p\leq n-p$, then $\Omega$ consists of all totally isotropic\footnote{In what follows we use the terms \texttt{isotropic} and \texttt{totally isotropic} as synonyms.} $k$-dimensional complex subspaces of $V$. See \cite{Wolf1} for more details in this latter case.

\subsection{The ind-group $\SL(\infty,\Cp)$ and its real forms}\fakesst In the rest of the paper, $V$ denotes a fixed countable-dimensional complex vector space with fixed basis $\Eo$. We fix an order on $\Eo$ via the ordered set $\Zp_{>0}$, i.e.,\break $\Eo=\{\epsi_1,~\epsi_2,~\ldots\}$. Let $V_*$ denote the span of the dual system $\Eo^*=\{\epsi_1^*,~\epsi_2^*,~\ldots\}$. By definition, the group $\GL(V,\Eo)$ is the group of invertible $\Cp$-linear transformations on $V$ that keep fixed all but finitely many elements of $\Eo$. It is not difficult to verify that that $\GL(V,\Eo)$ depends only on the pair $(V,V_*)$ but not on $\Eo$. Clearly, any operator from $\GL(V,\Eo)$ has a well-defined determinant. By $\SL(V,\Eo)$ we denote the subgroup of $\GL(V,\Eo)$ of all operators with determinant $1$. In the sequel $G = \SL(V,\Eo)$
and we also write $\SL(\infty,\Cp)$ instead of $G$.

Express the basis $\Eo$ as a union $\Eo=\bigcup \Eo_n$ of nested finite
subsets. Then $V$ is exhausted by the finite-dimensional subspaces
$V_n=\langle \Eo_n\rangle_{\Cp}$, i.e., $V=\varinjlim V_n$. To each linear operator $\vfi$ on~$V_n$ one can assign the operator $\wt\vfi$ on $V_{n+1}$ such that\break $\wt\vfi(x)=\vfi(x)\text{ for }x\in V_n,~\wt\vfi(\epsi_m)=\epsi_m\text{ for }\epsi_m\notin V_n$. This gives embeddings\break $\SL(V_n)\hookrightarrow\SL(V_{n+1})$, so that $G=\SL(V,\Eo)=\varinjlim\SL(V_n)$.
In what follows we consider this exhaustion of $G$ fixed, and set $G_n=\SL(V_n)$.

Recall that an \emph{ind-variety} over $\Rp$ or $\Cp$ (resp., an \emph{ind-manifold}) is an inductive limit of algebraic varieties (resp., of manifolds): $Y=\ilm Y_n$. Below we always assume that $Y_n$ form an ascending chain $$Y_1\hookrightarrow Y_2\hookrightarrow\ldots\hookrightarrow Y_n\hookrightarrow Y_{n+1}\hookrightarrow\ldots,$$ where $Y_n\hookrightarrow Y_{n+1}$ are closed embeddings. Any ind-variety or ind-manifold is endowed with a topology by declaring a subset $U\subset Y$ open if $U\cap Y_n$ is open for all $n$ in the corresponding topologies. A~\emph{morphism} $f\colon Y=\ilm Y_n\to Y'=\ilm Y_n'$ is a map induced by a collection of morphisms $\{f_n\colon Y_n\to Y_n'\}_{n\geq1}$ such that the restriction of $f_{n+1}$ to $Y_n$ coincides with $f_n$ for all $n\geq1$. A morphism $f\colon Y\to Y'$ is an \emph{isomorphism} if there exists a morphism $g\colon Y'\to Y$ for which $f\circ g=\id_{Y'}$ and $g\circ f=\id_Y$, where $\id$ is a morphism induced by the collection of the identity maps.

A \emph{locally linear algebraic ind-group} is an ind-variety $\Go=\bigcup\Go_n$ such that all $\Go_n$ are linear algebraic groups and the inclusions are group homomorphisms. In what follows we write \emph{ind-group} for brevity. Clearly, $G$ is an ind-group. By an \emph{ind-subgroup} of $G$ we understand a subgroup of $G$ closed in the direct limit Zariski topology. By definition, a real ind-subgroup $\Go^0$ of $G$ is called a \emph{real form} of~$G$, if $G$ can be represented as an increasing union $G=\bigcup\Go_n$ of its finite-dimensional Zariski closed subgroups such that $\Go_n$ is a semi-simple algebraic group and $\Go^0\cap\Go_n$ is a real form of~$\Go_n$ for each $n$. Below we recall the~classification of real forms of $G$ due to A. Baranov \cite{Baranov1}.

Fix a real structure $\tau$ on $V$ such that $\tau(e)=e$ for all $e\in\Eo$. Then each $V_n$ is $\tau$-invariant. Denote by $\GL(V_n,\Rp)$ (resp., by~$\SL(V_n,\Rp)$) the group of invertible (resp., of determinant 1) operators on $V_n$ defined over $\Rp$. Recall that a linear operator on a complex vector space with a real structure is \emph{defined over} $\Rp$ if it commutes with the real structure, or, equivalently, if it maps the real form to itself. For each $n$, the map $\vfi\mapsto\wt\vfi$ gives an embedding $\SL(V_n,\Rp)\hookrightarrow\SL(V_{n+1},\Rp)$, so the direct limit $G^0=\ilm\SL(V_n,\Rp)$ is well defined. We denote this real form of $G$ by $\SL(\infty,\Rp)$.

Fix a nondegenerate Hermitian form $\omega$ on $V$. Suppose that its restriction $\omega_n=\restr{\omega}{V_n}$ is nondegenerate for all $n$, and that $\omega(e_m,V_n)=0$ for $e_m\notin V_n$. Denote by $p_n$ the dimension of a maximal $\omega_n$-positive definite subspace of $V_n$, and put $q_n=\dim V_n-p_n$. Let $\SU(p_n,q_n)$ be the subgroup of $G_n$ consisting of all operators preserving the form $\omega_n$. For each $n$, the map $\vfi\mapsto\wt\vfi$ induces an embedding $\SU(p_n,q_n)\hookrightarrow\SU(p_{n+1},q_{n+1})$, so we have a direct limit $G^0=\ilm\SU(p_n,q_n)$. If there exists $p$ such that $p_n=p$ for all sufficiently large $n$ (resp., if $\lim_{n\to\infty}p_n=\lim_{n\to\infty}q_n=\infty$), then we denote this real form of $G$ by $\SU(p,\infty)$ (resp., by $\SU(\infty,\infty)$).

Finally, fix a \emph{quaternionic structure} $J$ on $V$, i.e., an antilinear automorphism of $V$ such that $J^2=-\id_V$. Assume that the complex dimension of $V_n$ is even for $n\geq1$, and that the restriction $J_n$ of $J$ to $V_n$ is a quaternionic structure on $V_n$. Furthermore, suppose that $$J(\epsi_{2i-1})=-\epsi_{2i},~J(\epsi_{2i})=\epsi_{2i-1}$$ for $i\geq1$. Let $\SL(V_n,\Hp)$ be the subgroup of $G_n$ consisting of all linear operators commuting with $J_n$, then, for each $n$, the map $\vfi\mapsto\wt\vfi$ induces an embedding of the groups $\SL(V_n,\Hp)\hookrightarrow\SL(V_{n+1},\Hp)$, and we denote the direct limit by $G^0=\SL(\infty,\Hp)=\ilm\SL(V_n,\Hp)$. This group is also a real form of $G$.

The next result is a corollary of \cite[Theorem 1.4]{Baranov1} and \cite[Corollary 3.2]{DimitrovPenkov2}.

\mtheo{If $G=\SL(\infty,\Cp)$\textup, then $\SL(\infty,\Rp)$\textup, $\SU(p,\infty)$\textup, $0\leq p<\infty$\textup, $\SU(\infty,\infty)$\textup, $\SL(\infty,\Hp)$ are all real\label{theo:list_of_real_forms} forms of $G$ up to isomorphism. These real forms are pairwise non-isomorphic as ind-groups.}

\subsection{Flag ind-varieties of the ind-group $G$}\fakesst \label{sst:gen_flags}Recall some basic definitions from \cite{DimitrovPenkov1}. A~\emph{chain} of subspaces in $V$ is a linearly ordered (by inclusion) set $\Co$ of distinct subspaces of $V$. We write $\Co'$ (resp., $\Co''$) for the subchain of $\Co$ of all $F\in\Co$ with an immediate successor (resp., an immediate predecessor). Also, we write $\Co^{\dag}$ for the set of all pairs $(F',~F'')$ such that $F''\in\Co''$ is the immediate successor of $F'\in\Co'$.

A \emph{generalized flag} is a chain $\Fo$ of subspaces in $V$ such that $\Fo=\Fo'\cup\Fo''$ and\break $V\setminus\{0\}=\bigcup\nolimits_{(F',F'')\in\Fo^{\dag}}F''\setminus F'$. Note that each nonzero vector $v\in V$ determines a unique pair $(F_v',F_v'')\in\Fo^{\dag}$ such that $v\in F''\setminus F'$. If $\Fo$ is a generalized flag, then each of $\Fo'$ and $\Fo''$ determines $\Fo$, because if $(F',F'')\in\Fo^{\dag}$, then
$F'=\bigcup\nolimits_{G''\in\Fo'',~G''\subsetneq F''}G''$, $F''=\bigcap\nolimits_{G'\in\Fo',~G'\supsetneq F'}G'$
(see \cite[Proposition 3.2]{DimitrovPenkov1}). We fix a linearly ordered set $(A,\preceq)$ and an isomorphism of ordered sets $A\to\Fo^{\dag}\colon a\mapsto(F_{\alpha}',F_{\alpha}'')$, so  that $\Fo$ can be written as $\Fo=\{F_{\alpha}',F_{\alpha}'',~\alpha\in A\}$. We will write $\alpha\prec\beta$ if $\alpha\preceq\beta$ and $\alpha\neq\beta$ for $\alpha,\beta\in A$.

A generalized flag $\Fo$ is called \emph{maximal} if it is not properly contained in another generalized flag. This is equivalent to the condition that $\dim F_v''/F_v'=1$ for all nonzero vectors $v\in V$. A generalized flag is called a \emph{flag} if the set of all proper subspaces of $\Fo$ is isomorphic as a linearly ordered set to a subset of $\Zp$.

We say that a generalized flag $\Fo$ is \emph{compatible} with a basis $E=\{e_1,~e_2,~\ldots\}$ of $V$ if there exists a surjective map $\sigma\colon E\to A$ such that every pair $(F_{\alpha}',F_{\alpha}'')\in\Fo^{\dag}$ has the form $F_{\alpha}'=\langle e\in E\mid \sigma(e)\prec\alpha\rangle_{\Cp}$, $F_{\alpha}''=\langle e\in E\mid \sigma(e)\preceq\alpha\rangle_{\Cp}$. By \cite[Proposition 4.1]{DimitrovPenkov1}, every generalized flag admits a compatible basis. A generalized flag $\Fo$ is \emph{weakly compatible} with $E$ if $\Fo$ is compatible with a basis $L$ of $V$ such that the set $E\setminus(E\cap L)$ is finite. Two generalized flags $\Fo$, $\Go$ are $E$-\emph{commensurable} if both of them are weakly compatible with $E$ and there exist an isomorphism of ordered sets $\phi\colon\Fo\to\Go$ and a finite-dimensional subspace $U\subset V$ such that
\begin{equation*}
\begin{split}
&\text{i) $\phi(F)+U=F+U$ for all $F\in\Fo$;}\\
&\text{ii) $\dim\phi(F)\cap U=\dim F\cap U$ for all $F\in\Fo$.}
\end{split}
\end{equation*}
Given a generalized flag $\Fo$ compatible with $E$, denote by $X=X_{\Fo,E}=\Fo\ell(\Fo,E)$ the set of all generalized flags in $V$, which are $E$-commensurable with $F$.

To endow $X$ with an ind-variety structure, fix an exhaustion $E=\bigcup E_n$ of $E$ by its finite subsets and denote $\Fo_n=\{F\cap \langle E_n\rangle_{\Cp},~F\in\Fo\}$. Given $\alpha\in A$, denote
\begin{equation*}
\begin{split}
d_{\alpha,n}'&=\dim F_{\alpha}'\cap\langle E_n\rangle_{\Cp}=|\{e\in E_n\mid\sigma(e)\prec\alpha\}|,\\
d_{\alpha,n}''&=\dim F_{\alpha}''\cap\langle E_n\rangle_{\Cp}=|\{e\in E_n\mid\sigma(e)\preceq\alpha\}|,
\end{split}
\end{equation*}
where $|\cdot|$ stands for cardinality. We define $X_n$ to be the projective varieties of flags in $\langle E_n\rangle_{\Cp}$ of the form $\{U_{\alpha}'',U_{\alpha}'',~\alpha\in A\}$, where $U_{\alpha}'$, $U_{\alpha}''$ are subspaces of $\langle E_n\rangle_{\Cp}$ of dimensions $d_{\alpha,n}'$, $d_{\alpha,n}''$ respectively, $U_{\alpha}'\subset U_{\alpha}''$ for all $\alpha\in A$, and $U_{\alpha}''\subset U_{\beta}'$ for all $\alpha\prec\beta$. (If $A$ is infinite, there exist infinitely many $\alpha$, $\beta\in A$ such that $U_{\alpha}=U_{\beta}$.) Define an embedding $\iota_n\colon X_n\to X_{n+1}\colon \{U_{\alpha}'',U_{\alpha}'',~\alpha\in A\}\mapsto\{W_{\alpha}'',W_{\alpha}'',~\alpha\in A\}$ by
\begin{equation}
\begin{split}\label{formula:iota_n}
W_{\alpha}'&=U_{\alpha}'\oplus\langle e\in E_{n+1}\setminus E_n\mid\sigma(e)\prec\alpha\rangle_{\Cp},\\
W_{\alpha}''&=U_{\alpha}''\oplus\langle e\in E_{n+1}\setminus E_n\mid\sigma(e)\preceq\alpha\rangle_{\Cp}.
\end{split}
\end{equation}
Then $\iota_n$ is a closed embedding of algebraic varieties, and there exists a bijection from $X$ to the inductive limit of this chain of morphisms, see \cite[Proposition 5.2]{DimitrovPenkov1} or \cite[Section 3.3]{FressPenkov1}. This bijection endows $X$ with an ind-variety structure which is independent on the chosen filtration $\bigcup E_n$ of the basis $E$. We will explain this bijection in more detail in Section~\ref{sect:orbits_as_manifolds}.

From now on we suppose that the linear span of $E_n$ coincides with $V_n$ and $V_*$ coincides with the span of the dual system $E^*=\{e_1^*,~e_2^*,~\ldots\}$. We assume also that the inclusion $G_n\hookrightarrow G_{n+1}$ induced by this exhaustion of $E$ coincides with the inclusion $\vfi\mapsto\wt\vfi$ defined above, i.e., that $\langle E_{n+1}\setminus E_n\rangle_{\Cp}=\langle\Eo_{n+1}\setminus\Eo_n\rangle_{\Cp}$. Denote by $H$ the ind-subgroup of $G=\SL(\infty,\Cp)$ of all operators from $G$ which are diagonal in $E$; $H$~is called a \emph{splitting Cartan subgroup} of $G$ (in fact, $H$~is a Cartan subgroup of $G$ in terminology of \cite{DimitrovPenkovWolf1}). We define a \emph{splitting Borel} (resp., \emph{parabolic}) \emph{subgroup} of $G$ to be and ind-subgroup of $G$ containing $H$ such that its intersection with $G_n$ is a Borel (resp., parabolic) subgroup of $G_n$. Note that if $P$ is a splitting parabolic subgroup of $G$ and $P_n=P\cap G_n$, then $G/P=\bigcup G_n/P_n$ is a locally projective ind-variety, i.e., an ind-variety exhausted by projective varieties. One can easily check that the group $G$ naturally acts on $X$. Given a generalized flag $\Fo$ in $V$ which is compatible with~$E$, denote by $P_{\Fo}$ the stabilizer of $\Fo$ in $G$. For the proof of the following theorem, see \cite[Proposition 6.1, Theorem 6.2]{DimitrovPenkov1}.
\mtheo{Let $\Fo$ be a generalized flag compatible with $E$, $X=\Fl(\Fo,E)$ and $G=\SL(\infty,\Cp)$.\\ \textup{i)}The group $P_{\Fo}$ is a parabolic subgroup of $G$ containing $H$\textup, and the map $\Fo\mapsto P_{\Fo}$ is a bijection between generalized flags compatible with $E$ and splitting parabolic subgroups of $G$.\\ \textup{ii)} The ind-variety $X$ is in fact $G$-homogeneous\textup, and the map $g\mapsto g\cdot\Fo$ induces an isomorphism of ind-varieties $G/P_{\Fo}\cong X$.\\ \textup{iii)} $\Fo$ is maximal if and only if $P_{\Fo}$ is a splitting Borel subgroup of $G$.}

\exam{i) A first example \label{exam:gen_flags}of generalized flags is provided by the flag $\Fo=\{\{0\}\subset F\subset V\}$, where $F$ is a proper nonzero subspace of $V$. If $F$ is compatible with $E$, then we can assume that $F=\langle\sigma\rangle_{\Cp}$ for some subset $\sigma$ of $E$. In this case the ind-variety $X$ is called an \emph{ind-grassmannian}, and is denoted by~$\Grr{F}{E}$. If $k=\dim F$ is finite, then a flag $\{\{0\}\subset F'\subset V\}$ is $E$-commensurable with $\Fo$ if and only if $\dim F=k$, hence $\Grr{F}{E}$ depends only on $k$, and we denote it by $\Grr{k}{V}$. Similarly, if $k=\codim_VF$ is finite, then $\Grr{F}{E}$ depends only on $E$ and $k$ (but not on~$F$) and is isomorphic to $\Grr{k}{V_*}$: an isomorphism $\Grr{F}{E}\to\{F\subset V_*\mid\dim F=k\}=\Gra{k,V_*}$ is induced by the map $\Grr{F}{E}\ni U\mapsto U^{\#}=\{\phi\in V_*\mid\phi(x)=0\text{ for all }x\in U\}$. Finally, if $F$ is both infinite dimensional and infinite codimensional, then $\Grr{F}{E}$ depends on $F$ and~$E$, but all such ind-varieties are isomorphic and denoted by $\Gra{\infty}$, see \cite{PenkovTikhomirov1} or \cite[Section 4.5]{FressPenkov1} for the details. Clearly, in each case one has $\Fo'=\{\{0\}\subset F\}$, $\Fo''=\{F\subset V\}$.

ii) Our second example is the generalized flag $\Fo=\{\{0\}=F_0\subset F_1\subset\ldots\}$, where $F_i=\langle e_1,\ldots,e_i\rangle_{\Cp}$ for all $i\geq1$. This clearly is a flag. A flag $\wt\Fo=\{\{0\}=\wt F_0\subset\wt F_1'\subset \ldots\}$ is $E$-commensurable with $\Fo$ if and only if $\dim F_i=\dim\wt F_i$ for all~$i$, and $F_i=\wt F_i$ for large enough $i$. The flag $\Fo$ is maximal, and $\Fo'=\Fo$, $\Fo''=\Fo\setminus\{0\}$.

iii) Put $\Fo=\{\{0\}=F_0\subset F_1\subset F_2\subset\ldots\subset F_{-2}\subset F_{-1}\subset V\}$, where\break $F_i=\langle e_1,~e_3,~\ldots, e_{2i-1}\rangle_{\Cp}$, $F_{-i}=\langle\{e_j,~j\text{ odd}\}\cup\{e_{2j},~j>i\}\rangle_{\Cp}$ for $i\geq1$. This generalized flag is clearly not a flag, and is maximal. Here $\Fo'=\Fo\setminus V$, $\Fo''=\Fo\setminus\{0\}$. Note also that $\wt\Fo\in X=\Fl(\Fo,E)$ does not imply that $\wt F_i=F_i$ for $i$ large enough. For example, let $\wt F_1=\Cp e_2$, $$\wt F_i=\langle e_2,~e_3,~e_5,~e_7,~\ldots,~e_{2i-1}\rangle_{\Cp}$$ for $i>1$, and $\wt F_{-i}=\langle \{e_j,~j\text{ odd},~j\geq3\}\cup\{e_2\}\cup\{e_{2j},j>i\}\rangle_{\Cp},~i\geq1$, then $\wt\Fo\in X$, but $\wt F_i\neq F_i$ for all $i$.}

\nota{In all above examples $X=G/P_{\Fo}$, where $P_{\Fo}$ is the stabilizer of $\Fo$ in $G$. The ind-grassmannians in (i) are precisely the ind-varieties $G/P_{\Fo}$ for maximal splitting parabolic ind-subgroups $P_{\Fo}\subset G$. The ind-variety $\Fl(\Fo,E)$, where $\Fo$ is the flag in (ii), equals $G/P_{\Fo}$ where $P_{\Fo}$ is the upper-triangular Borel ind-subgroup in the realization of $G$ as $\Zp_{>0}\times\Zp_{>0}$-matrices.}

\section{$G^0$-orbits as\label{sect:orbits_as_manifolds} ind-manifolds}\fakesect

In this section, we establish a basic property of the orbits on $G/P$ of a real form $G^0$ of\break $G=\SL(\infty,\Cp)$. Precisely, we prove that the intersection of a $G^0$-orbit with $X_n$ is a single orbit. Consequently, each $G^0$-orbit is an infinite-dimensional real ind-manifold.

We start by describing explicitly the bijection $X\to\ilm X_n$ mentioned in Subsection~\ref{sst:gen_flags}. Let $\Fo$ be a generalized flag in $V$ compatible with the basis $E$, and $X=\Fl(\Fo,E)$ be the corresponding ind-variety of generalized flags. Recall that we consider $X$ as the inductive limit of flag varieties $X_n$, where the embeddings $\iota_n\colon X_n\hookrightarrow X_{n+1}$ are defined in the previous subsection. Put $E_m'=\{e_1,~e_2,~\ldots,~e_m\}$ and $\Vo_m=\langle E_m'\rangle_{\Cp}$. The construction of $\iota_n$ can be reformulated as follows.

The dimensions of the spaces of the flag $\Fo\cap\Vo_m$ form a sequence of integers $$0=d_{m,0}<d_{m,1}<\ldots<d_{m,s_{m-1}}<d_{s_m}=\dim\Vo_m=m.$$ Let $\Fl(d_m,\Vo_m)$ be the flag variety of type $d_m=(d_{m,1},\ldots,d_{m,s_{m-1}})$ in $\Vo_m$. Since either $s_{m+1}=s_m$ or $s_{m+1}=s_m+1$, there is a unique $j_m$ such that $d_{m+1,i}=d_{m,i}+1$ for $0\leq i<j_m$ and $d_{m+1,j_m}>d_{m,j_m}$. Then, for $j_m\leq i<s_m$, $d_{m+1,i}=d_{m,i}+1$ in case $s_{m+1}=s_m$, and $d_{m+1,i}=d_{m,i-1}+1$ in case $s_{m+1}=s_m+1$. In other words, $j_m\leq s_m$ is the minimal nonnegative integer for which there is $\alpha\in A$ with $$\dim F_{\alpha}''\cap\Vo_{m+1}=\dim F_{\alpha}''\cap\Vo_m+1.$$

Now, for each $m$ we define an embedding $\xi_m\colon\Fl(d_m,\Vo_m)\hookrightarrow\Fl(d_{m+1},\Vo_{m+1})$: given a flag $\Go_m=\{\{0\}=G_0^m\subset G_1^m\subset\ldots\subset G_{s_m}^m=V_m\}\in\Fl(d_m,\Vo_m)$, we set $\xi_m(\Go_m)=\Go_{m+1}=\{\{0\}=G_0^{m+1}\subset G_1^{m+1}\subset\ldots\subset G_{s_{m+1}}^{m+1}=V_{m+1}\}\in\Fl(d_{m+1},\Vo_{m+1})$, where
\begin{equation}
G_i^{m+1}=\begin{cases}
G_i^m,&\text{if\label{formula:iota_n_another} $0\leq i<j_m$},\\
G_i^m\oplus\Cp e_{m+1},&\text{if $j_m\leq i\leq s_{m+1}$ and $s_{m+1}=s_m$},\\
G_{i-1}^m\oplus\Cp e_{m+1},&\text{if $j_m\leq i\leq s_{m+1}$ and $s_{m+1}=s_m+1$}.
\end{cases}
\end{equation}

For any $\Go\in X$ we choose a positive integer~$m_{\Go}$ such that $\Fo$ and $\Go$ are compatible with bases containing $\{e_i\mid ~i\geq m_{\Go}\}$, and $\Vo_{m_{\Go}}$ contains a subspace which makes these generalized flags $E$-commensurable. In addition, we can assume that\break $m_{\Fo}\leq m_{\Go}$ for all $\Go\in X$ (in fact, we can set $m_{\Fo}=1$ because $\Fo$ is compatible with $E$). Let $m_{\Fo}\leq m_1<m_2<\ldots$ be an arbitrary sequence of integer numbers. For $n\geq1$, denote $E_n=E_{m_n}'$, $V_n=\Vo_{m_n}$. Then $X_n=\Fl(d_{m_n},\Vo_{m_n})$ and, according to (\ref{formula:iota_n}), $\iota_n=\xi_{m_{n+1}-1}\circ\xi_{m_{n+1}-2}\circ\ldots\circ\xi_{m_n}$. The bijection $X\to\ilm X_n$ from Subsection~\ref{sst:gen_flags} now has the form $\Go\mapsto\ilm\Go_n$, where $\Go_n=\{F\cap V_n,~F\in\Go\}$ for $n$ such that $m_n\geq m_{\Go}$. By a slight abuse of notation, in the sequel we will denote the canonical embedding $X_n\hookrightarrow X$ by the same letter $\iota_n$.

Let $G^0$ be a real form of~$G=\SL(\infty,\Cp)$ (see Theorem \ref{theo:list_of_real_forms}). The group\break $G_n=\SL(V_n)$ naturally acts on $X_n$, and the map $\iota_n$ is equivariant: $g\cdot\iota_n(x)=\iota_n(g\cdot x)$, $g\in G_n\subset G_{n+1}$, $x\in X_n$. Put also $G_n^0=G^0\cap G_n$. Then $G_n^0$ is a real form of~$G_n$. For the rest of the paper we fix some specific assumptions on $V_n$ for different real forms. We now describe these assumptions case by case.

Let $G^0=\SU(p,\infty)$ or $\SU(\infty,\infty)$. Recall that the restriction $\omega_n$ of the fixed nondegenerate Hermitian form $\omega$ to $V_n$ is nondegenerate. From now on, we assume that if $e\in E_{n+1}\setminus E_n$, then $e$ is orthogonal to~$V_n$ with respect to $\omega_{n+1}$. Next, let $G^0=\SL(\infty,\Rp)$. Here we assume that $m_n$ is odd for each $n\geq1$, and that $\langle E_n\rangle_{\Rp}$ is a real form of $V_n$. Finally, for $G^0=\SL(\infty,\Hp)$, we assume that $m_n$ is even for all $n\geq1$ and that $J(e_{2i-1})=-e_{2i}$, $J(e_{2i})=e_{2i-1}$ for all $i$. These additional assumptions align the real form $G^0$ with the flag variety $X$.

Our main result in this section is as follows.

\theop{If \label{theo:Omega_n_cap}
$\iota_n(X_n)$ has nonempty intersection with a $G_{n+1}^0$-orbit,
then that intersection is a single $G_n^0$-orbit.}
{The proof goes case by case.

\textsc{\underline{Case $G^0=\SU(\infty,\infty)$}}. (The proof for $G^0=\SU(p,\infty)$, $0\leq p<\infty$, is completely similar.) Pick two flags
\begin{equation*}
\begin{split}
\Au&=\{\{0\}=A_0\subset A_1\subset\ldots\subset A_{s_{m_n}}=V_n\},\\
\Bu&=\{\{0\}=B_0\subset B_1\subset\ldots\subset B_{s_{m_n}}=V_n\}
\end{split}
\end{equation*}
in $X_n$ such that $\wt\Au=\iota_n(\Au)$ and $\wt\Bu=\iota_n(\Bu)$ belong to a given $G_{n+1}^0$-orbit.

Put
\begin{equation*}\predisplaypenalty=10000
\begin{split}
\wt\Au&=\{\{0\}=\wt A_0\subset \wt A_1\subset\ldots\subset\wt A_{s_{m_{n+1}}}=V_{n+1}\},\\
\wt\Bu&=\{\{0\}=\wt B_0\subset \wt B_1\subset\ldots\subset\wt B_{s_{m_{n+1}}}=V_{n+1}\}.
\end{split}
\end{equation*}
There exists $\vfi\in\SU(\omega_{n+1},V_{n+1})$ satisfying $\vfi(\wt\Au)=\wt\Bu$, i.e., $\vfi(\wt A_i)=\wt B_i$ for all $i$ from $0$ to~$s_{m_{n+1}}$. To prove the result, we must construct an isometry $\vfi\colon V_n\to V_n$ satisfying $\vfi(\Au)=\Bu$. Of course, one can scale $\vfi$ to obtain an isometry of determinant $1$. By Huang's extension of Witt's Theorem \cite[Theorem 6.2]{Huang1}, such an isometry exists if and only if $A_i$ and $B_i$ are isometric for all $i$ from $1$ to $s_{m_n}$, and
\begin{equation}
\dim(A_i\cap A_j^{\perp,V_n})=\dim(B_i\cap B_j^{\perp,V_n})\label{formula:Witt_for_flags}
\end{equation}
for all $i<j$ from $1$ to $s_{m_n}$. (Here $U^{\perp,V_n}$ denotes the $\omega_n$-orthogonal complement within $V_n$ of a subspace $U\subset V_n$.) Pick $i$ from $1$ to $s_{m_n}$. Since $e_{n+1}$ is orthogonal to $V_n$ and $\wt\vfi$ establishes an isometry between $\wt A_i$ and $\wt B_i$, the first condition is satisfied. So it remains to prove (\ref{formula:Witt_for_flags}).

To do this, denote $C_n=\langle E_{n+1}\setminus E_n\rangle_{\Cp}$. Since $C_n$ is orthogonal to $V_n$, for given subspaces $U\subset V_n$, $W\subset C_n$ one has $(U\oplus W)^{\perp,V_{n+1}}=U^{\perp,V_n}\oplus W^{\perp,C_n}$. Hence, if\break $\wt A_k=A_k\oplus W_k$, $\wt B_k=B_k\oplus W_k$ for $k\in\{i,j\}$ and some subspaces of $W_i$, $W_j\subset C_n$, then $$\wt A_i\cap\wt A_j^{\perp,V_{n+1}}=(A_i\oplus W_i)\cap(A_j^{\perp,V_n}\oplus W_j^{\perp,C_n})=(A_i\cap A_j^{\perp,V_n})\oplus(W_i\cap W_j^{\perp,C_n}),$$ and the similar equality holds for $\wt B_i\cap\wt B_j^{\perp,V_{n+1}}$. The result follows.

\textsc{\underline{Case $G^0=\SL(\infty,\Rp)$}}. Here we first prove that if $\Au$ and $\Bu$ are flags in $\Vo_n$,\break $\wt\Au$ and $\wt\Bu$ are their images in $\Vo_{n+1}$ under the map $\xi_n$, and there exists\break $\vfi\in\GL(\Vo_{n+1},\Rp)$ satisfying $\vfi(\wt\Au)=\wt\Bu$, then there exists an operator $\nu\in\GL(\Vo_n,\Rp)$ such that $\nu(\Au)=\Bu$.

Consider first the case when $\vfi(e_{n+1})\notin\Vo_n$. Denote $\vfi(e_{n+1})=v+te_{n+1}$, $v\in\Vo_n$, $t\in\Rp$, $t\neq0$. Then $t^{-1}\vfi\in\GL(\Vo_{n+1},\Rp)$ maps $\wt\Au$ to $\wt\Bu$, so we can assume that $t=1$, i.e., $\vfi(e_{n+1})=v+e_{n+1}$. Since $$\vfi(A_{j_n}\oplus\Cp e_{n+1})=\vfi(\wt A_{j_n})=\wt B_{j_n}=B_{j_n}\oplus\Cp e_{n+1},$$ the vector $v$ belongs to $B_i$ for all $i\geq j_n$. Let $\psi\in\GL(\Vo_{n+1},\Rp)$ be defined by $$\psi(x+se_{n+1})=x+s(e_{n+1}-v),~x\in\Vo_n,~s\in\Cp.$$ Clearly, $\psi(\vfi(e_{n+1}))=e_{n+1}$.

If $i<j_n$ and $x\in A_i$, then $\vfi(x)\in B_i\subset\Vo_n$, so $\psi(\vfi(x))=\vfi(x)\in B_i$. If $i\geq j_n$ and $x\in A_r$, where $r=i$ for $s_{n+1}=s_n$ and $r=i-1$ for $s_{n+1}=s_n+1$, then we put $\vfi(x)=y+se_{n+1}$, $y\in B_i$, $s\in\Cp$. One has $$\psi(\vfi(x))=\psi(y+se_{n+1})=y+s(e_{n+1}-v)\in B_r\oplus\Cp e_{n+1}=\wt B_i.$$ In both cases the operator $\psi\circ\vfi$ maps
$\wt A_i$ to $\wt B_i$ for all $i$ from $0$ to  $s_{n+1}$. Hence we may assume without loss of generality that $\vfi(e_{n+1})=e_{n+1}$. Then the operator\break $\nu=\pi\circ\restr{\vfi}{\Vo_n}$, where $\pi\colon\Vo_{n+1}\to\Vo_n$ is the projection onto $\Vo_n$ along $\Cp e_{n+1}$, is invertible, is defined over $\Rp$, and maps each $A_i$ to $B_i$, $0\leq i\leq s_n$, as required.

Suppose now that $\vfi(e_{n+1})=b\in\Vo_n$. In this case $s_{n+1}=s_n$ because the condition $$\vfi(A_{j_n-1}\oplus\Cp e_{n+1})=\vfi(\wt A_{j_n})=\wt B_{j_n}=B_{j_n-1}\oplus\Cp e_{n+1}$$ contradicts the equality $s_{n+1}=s_n+1$. Arguing as above, we see that $b\in B_i$ for all $i\geq j_n$. If $\vfi^{-1}(e_{n+1})=a\notin\Vo_n$, then one can construct $\nu$ as in the case when $b\notin\Vo_n$ with $\vfi^{-1}$ instead of~$\vfi$. Therefore, we may assume that $a\in A_i$ for all $i\geq j_n$. Let $U^0$ be an $\Rp$-subspace of $\Vo_n^0$ such that $\Vo_n^0=U^0\oplus\Rp b$, then $\Vo_n=U\oplus\Cp b$, where\break $U=\Cp\otimes_{\Rp}U^0$. If $a,~b$ are linearly independent, we choose $U^0$ so that $a\in U^0$. Define $\nu$ as follows: if $\vfi(x)=y+sb+re_{n+1}$, $x\in\Vo_n$, $y\in U$, $s,r\in\Cp$, then put\break $\nu(x)=y+(s+r)b$. One can easily check that $\nu$ satisfies all required conditions.

Now we are ready to prove the result for $G^0=\SL(\infty,\Rp)$. Namely, let $\Au$, $\Bu\in X_n$, and $\vfi\in\SL(V_n,\Rp)$ satisfy $\vfi(\iota_n(\Au))=\iota_n(\Bu)$, then $\vfi$ belongs to $\GL(\Vo_{m_{n+1}},\Rp)$. Hence there exists $\nu'\in\GL(\Vo_{m_{n+1}-1},\Rp)$ which maps $\xi_{m_{n+1}-2}\circ\ldots\circ\xi_{m_n}(\Au)$ to $\xi_{m_{n+1}-2}\circ\ldots\circ\xi_{m_n}(\Bu)$ because $\iota_n=\xi_{m_{n+1}-1}\circ\xi_{m_{n+1}-2}\circ\ldots\circ\xi_{m_n}$. Continuing this process, we see that there exists an operator $\nu''\in\GL(V_n,\Rp)$ such that\break $\nu''(\Au)=\Bu$. Since $V_n$ is odd-dimensional, one can scale $\nu''$ to obtain a required operator $\nu\in\SL(V_n,\Rp)$.

\textsc{\underline{Case $G^0=\SL(\infty,\Hp)$}}. Let $\Au$, $\Bu$ be two flags in $\Vo_{2n}$ and $\wt\Au=\xi_{2n+1}\circ\xi_{2n}(\Au)$, $\wt\Bu=\xi_{2n+1}\circ\xi_{2n} (\Bu)$. Let $\vfi\in\SL(\Vo_{2n+2},\Hp)$ satisfy $\vfi(\wt\Au)=\wt\Bu$. Our goal is to construct $\nu\in\SL(\Vo_{2n},\Hp)$ such that $\nu(\Au)=\Bu$. Then, repeated application of this procedure will imply the result.

For simplicity, denote $e=e_{2n+1}$, $e'=e_{2n+2}$. Recall that $J(e)=-e'$, $J(e')=e$, and note that $b=\vfi(e)\in\Vo_{2n}$ if and only if $b'=\vfi(e')\in\Vo_{2n}$, because $\Vo_{2n}$ is $J$-invariant and $\vfi$ commutes with $J$.

First, suppose that both $b$~and ~$b'$ do not belong to $\Vo_{2n}$. The vector $b$ admits a unique representation in the form $b=v+te+t'e'$ for $v\in\Vo_{2n}$, $t,~t'\in\Cp$. Then $$b'=\vfi(e')=\vfi(-J(e))=-J(\vfi(e))=-J(b)=v'-\bar t'e+\bar te',$$ where $v'=-J(v)\in\Vo_{2n}$. Set $$T=\begin{pmatrix}t&-\bar t'\\t'&\bar t\end{pmatrix},~d=\det T=|t|^2+|t'|^2\in\Rp_{>0}.$$ Let $\psi\in\GL(\Vo_{2n+2})$ be the operator defined by $\psi(x)=x$, $x\in\Vo_{2n}$,
\begin{equation*}
\begin{split}
&\psi(e)=-d^{-1}(\bar t(v+e)-t'(v'+e')),\\
&\psi(e')=-d^{-1}(\bar t'(v+e)+t(v'+e')).
\end{split}
\end{equation*}
It is easy to see that $\psi$ commutes with $J$, $\det\psi=\det T^{-1}\in\Rp_{>0}$, and $\psi(b)=e$, $\psi(b')=e'$. Furthermore, one can check that $\nu=\pi\circ\psi\circ\restr{\vfi}{\Vo_{2n}}\colon\Vo_{2n}\to\Vo_{2n}$ commutes with $J$ and maps $\Au$ to $\Bu$, where $\pi\colon V_{n+1}\to V_n$ is the projection onto $V_n$ along $\Cp e\oplus\Cp e'$. Since $\det\nu\in\Rp_{>0}$, one can scale $\nu$ to obtain an operator from $\SL(\Vo_{2n},\Hp)$, as required.

Second, suppose that $b, b'\in\Vo_{2n}$. If $a=\vfi^{-1}(e)$ and $a'=\vfi^{-1}(e')$ do not belong to $V_n$, one can argue as in the first case with $\vfi^{-1}$ instead of $\vfi$, so we may assume without loss of generality that $a,a'\in V_n$. (Note that if $a,~a',~b,~b'$ are linearly dependent, then $\Cp a\oplus\Cp a'=\Cp b\oplus\Cp b'$.) In this case, denote by $U$ a $J$-invariant subspace of $V_n$ spanned by some basic vectors $e_i$ such that $V_n=U\oplus\Cp b\oplus\Cp b'$. (If $a,~a',~b,~b'$ are linearly independent, we choose $U$ such that $a,a'\in U$.) Define~$\nu$ by the following rule: if $\vfi(x)=y+sb+sb'+re+r'e'$, $x\in V_n$, $y\in U$, $s,~s',~r,~r'\in\Cp$, then $\nu(x)=y+(s+r)b+(s'+r')b'$. One can check that $\det\nu=\det\vfi=1$, $\nu$ commutes with $J$ (so $\nu\in\SL(V_n,\Hp)$) and maps each $A_i$, $0\leq i\leq s_n$, to $B_i$. Thus, $\nu$ satisfies all required conditions.}

The following result is an immediate corollary of this theorem.

\corop{Let $\Omega$ be a $G^0$-orbit on $X$\textup, and $\Omega_n=\iota_n^{-1}(\Omega)\subset X_n$. Then\\\textup{i)} $\Omega_n$ is a single $G_n^0$-orbit\textup;\\\textup{ii)} $\Omega$ is an infinite-dimensional real ind-manifold.}{i) Suppose $\Au$, $\Bu\in\Omega_n$. Then there exists $m\geq n$ such that images of $\Au$ and $\Bu$ under the morphism $\iota_{m-1}\circ\iota_{m-2}\circ\ldots\circ\iota_n$ belong to the same $G_m^0$-orbit. Applying Theorem~\ref{theo:Omega_n_cap} subsequently to $\iota_{m-1}$, $\iota_{m-2}$, $\ldots$, $\iota_n$, we see that $\Au$ and~$\Bu$ belong to the same $G_n^0$-orbit.

ii) By definition, $\Omega=\ilm\Omega_n$. Next, (i) implies that $\Omega$ is a real ind-manifold. By Theorem~\ref{theo:finite_dim_case}~(v), we have $\dim_{\Rp}\Omega_n\geq\dim_{\Cp}X_n$. Since $\lim_{n\to\infty}\dim_{\Cp}X_n=\infty$, we conclude that $\Omega$ is infinite dimensional.}

\section{Case of finitely many $G^0$-orbits}\fakesect

We give now a criterion for $X=\Fl(\Fo,E)$ to have a finite number of $G^0$-orbits, and observe that, if this is the case, the degeneracy order on the $G^0$-orbits in $X$ coincides with that on the $G_n^0$-orbits in $X_n$ for large enough $n$. Recall that the \emph{degeneracy order} on the orbits is the partial order $\Omega\leq\Omega'\iff\Omega\subseteq\overline{\Omega}'$.

A generalized flag $\Fo$ is \emph{finite} if it consists of finitely many (possibly infinite-dimensional) subspaces. We say that a generalized flag $\Fo$ has \emph{finite type} if it consists of finitely many subspaces of $V$ each of which has either finite dimension or finite codimension in $V$. A finite type generalized flag is clearly a flag. An ind-variety $X=\Fl(\Fo,E)$ is \emph{of finite type} if $\Fo$ is of finite type (equivalently, if any $\wt\Fo\in X$ is of finite type).

\propp{For $G^0=\SU(\infty,\infty)$, $\SL(\infty,\Rp)$ and $\SL(\infty,\Hp)$, the number of $G^0$-orbits\label{prop:finiteness} on $X$ is finite if and only if $X$ is of finite type. For $G^0=\SU(p,\infty)$, $0<p<\infty$, the number of $G^0$-orbits on $X$ is finite if and only if $\Fo$ is finite. For $G^0=\SU(0,\infty)$, the number of $G^0$-orbits on $X$ equals $1$.}
{

\textsc{\underline{Case $G^0=\SU(\infty,\infty)$}}. First consider the case $X=\Grr{F}{E}$, where $F$ is a subspace of $V$. Clearly, $X$ is of finite type if and only if $\dim F<\infty$ or $\codim_VF<\infty$. Note that for ind-grassmannians, the construction of $\iota_n$ from (\ref{formula:iota_n}) is simply the following. Given $n$, let $W_{n+1}$ be the span of $E_{n+1}\setminus E_n$, and $U_{n+1}$ be a fixed $(k_{n+1}-k_n)$-dimensional subspace of $W_{n+1}$, where $k_i=\dim F\cap V_i$. Then the embedding\break $\iota_n\colon X_n=\Grr{k_n}{V_n}\to X_{n+1}=\Grr{k_{n+1}}{V_{n+1}}$ has the form $\iota_n(A)=A\oplus U_{n+1}$ for $A\in X_n$.

Recall that if $\codim_VF=k$, then the map $$U\mapsto U^{\#}=\{\phi\in V_*\mid\phi(x)=0\text{ for all }x\in U\}$$ induces an isomorphism $\Grr{F}{E}\to\{F'\subset V_*\mid\dim F'=k\}=\Grr{k}{V_*}$; we denote this isomorphism by $D$. To each operator $\psi\in\GL(V,E)$ one can assign the linear operator $\psi_*$ on $V_*$ acting by $(\psi_*(\lambda))(x)=\lambda(\psi(x))$, ${\lambda\in V_*}$, $x\in V$. This defines an isomorphism $\SL(V,E)\to\SL(V_*,E^*)$, and $D$ becomes a $G$-equivariant isomorphism of ind-varieties. Hence, for $X$ of finite type, we can consider only the case when $\dim F=k$.

If $\dim F=k$, then $X$ consists of all $k$-dimensional subspaces of $V$. Pick $A$, $B\in X$. There exists $n$ such that $X_n=\Grr{k}{V_n}$ and $A,~B\in\iota_n(X_n)$. Witt's Theorem shows that, for each $m\geq n$, $A$ and $B$ belong to the same $G_m^0$-orbit if and only if the signatures of the forms $\restr{\omega_m}{A}$ and $\restr{\omega_m}{B}$ coincide. Since $\restr{\omega_m}{A,B}=\restr{\omega}{A,B}$, we conclude $A$ and $B$ belong to the same $G^0$-orbit if and only if their signatures coincide. Thus, the number of $G^0$-orbits on $X$ is finite.

On the other hand, if $\dim F=\codim_VF=\infty$, then $$\lim_{n\to\infty}k_n=\lim_{n\to\infty}(\dim V_n-k_n)=\infty.$$ In this case, the number of possible signatures of the restriction of $\omega_n$ to a $k_n$-dimensional subspace tends to infinity, hence the number of $G_n^0$-orbits tends to infinity. By Theorem~\ref{theo:Omega_n_cap}, the number of $G^0$-orbits on~$X$ is infinite.

Now, consider the general case $X=\Fl(\Fo,E)$. Let $\Fo$ be of finite type. Then\break $\Fo=\Au\cup\Bu$ where $\Au$~and $\Bu$ are finite type subflags of $\Fo$ consisting of finite-dimensional and finite-codimensional subspaces from $\Fo$ respectively. Note that $\Au$ and $\Bu$ are compatible with the basis $E$, hence there exists $N$ such that if $n\geq N$, then $A\subseteq V_n$ for all $A\in\Au$ and $\codim_{V_n}(B\cap V_n)=\codim_VB$ for all $B\in\Bu$. Set
\begin{equation*}
\begin{split}
\Au&=\{A_1\subset A_2\subset\ldots\subset A_k\},\\
\Bu&=\{B_1\subset B_2\subset\ldots\subset B_l\},
\end{split}
\end{equation*}
and $a_i=\dim A_i$, $1\leq i\leq k$, $b_i=\codim_VB_i$, $1\leq i\leq l$.

Denote by $s(U)$ the signature of $\restr{\omega}{U}$ for a finite-dimensional subspace $U\subset V$. According to~\cite[Theorem 6.2]{Huang1}, to check that the number of $G^0$-orbits on $X$ is finite, it is enough to prove that all of the following sets are finite:
\begin{equation*}
\begin{split}
S_A&=\{s(A)\mid A\subset V_n,~n\geq N,~\dim A=a_i\text{ for some }i\},\\
S_B&=\{s(B)\mid B\subset V_n,~n\geq N.~\codim_{V_n}B=b_i\text{ for some }i\},\\
P_A&=\{\dim A\cap A_0^{\perp,V_n}\mid A,~A_0\subset V_n,~n\geq N,\\
&\dim A=a_i,~\dim A_0=a_j\text{ for some }i<j\},\\
P_B&=\{\dim B\cap B_0^{\perp,V_n}\mid B,~B_0\subset V_n,~n\geq N,\\
&\codim_{V_n}B=b_i,~\codim_{V_n}B_0=b_j\text{ for some }i<j\},\\
P_{AB}&=\{\dim A\cap B^{\perp,V_n}\mid A,~B\subset V_n,~n\geq N,~\\
&\dim A=a_i,~\codim_{V_n}B=b_j\text{ for some }i,j\}.
\end{split}
\end{equation*}
The finiteness of $S_A$ and $P_A$ is obvious. In particular, this implies that the number of $G^0$-orbits on $\Fl(\Au,E)$ is finite. Applying the map  $U\mapsto U^{\#}$ described above, we see that the number of $G^0$-orbits on $\Fl(\Bu,E)$ is finite. Consequently, the sets $S_B$ and $P_B$ are finite. Finally, since $\omega_n=\restr{\omega}{V_n}$ is nondegenerate for each $n$, we see that if $B\subset V_n$ and $\codim_{V_n}B=b_i$ for some $i$, then $\dim B^{\perp,V_n}=\codim_{V_n}B=b_i$. Hence $P_{AB}$ is finite. Thus, if $\Fo$ is of finite type then
the number of $G^0$-orbits on $\Fl(\Fo,E)$ is finite.

On the other hand, suppose that $\Fo$ is not of finite type. If there is a space $F\in\Fo$ with $\dim F=\codim_VF=\infty$, then we are done, because the map $$X\to\Grr{F}{E}\colon \Go\mapsto\text{the subspace in $\Go$ corresponding to $F$}$$ is a $G$-equivariant epimorphism of ind-varieties, and the number of $G^0$-orbits on the ind-grassmannian $\Grr{F}{E}$ is infinite by the above.

If all $F\in\Fo$ are of finite dimension or finite codimension, there exist subspaces $F_n\in\Fo$ of arbitrarily large dimension or arbitrarily large codimension. In the former case the statement follows from the fact that the number of possible signatures of such spaces tends to infinity, and in the latter case the statement gets reduced to the former one via the map $U\mapsto U^{\#}$.

\textsc{\underline{Case $G^0=\SU(p,\infty)$, $0<p<\infty$.}} First suppose that $\Fo$ is finite, i.e., $|\Fo|=N<\infty$. Given $n\geq1$, denote $S_n=\{s(A)\mid A\subset V_n\}$ and $P_n=\{\dim A\cap B^{\perp,V_n}\mid A\subset B\subset V_n\}$. Let $s(A)=(a,b,c)$ for some subspace $A$ of $V_n$. Then, clearly, $a\leq p$ and $c\leq p$, hence $|S_n|\leq p^2$. On the other hand, if $A\subset B$ are subspaces of $V_n$ then $A^{\perp,V_n}\supset B^{\perp,V_n}$, so $A\cap B^{\perp,V_n}\subset A\cap A^{\perp,V_n}$. But $\dim A\cap A^{\perp,V_n}=c\leq p$. Thus $|P|\leq p$. Now \cite[Theo\-rem~6.2]{Huang1} shows that the number of $G_n^0$-orbits on $X_n$ is less or equal to ${N|S_n|N^2|P_n|\leq N^3p^3}$. Hence, by Theorem~\ref{theo:Omega_n_cap}, the number of $G^0$-orbits on $X$ is finite.

Now suppose that $\Fo$ is infinite. In this case, given $m\geq1$, there exists $n$ such that the length of each flag from $X_n$ is not less than~$m$, the positive index of $\restr{\omega}{V_n}$ (i.e., the dimension of a maximal positive definite subspace of $V_n$) equals $p$, and $\codim_{V_n}F_m\geq p$, where $\Fo_n=\{F_1\subset\ldots\subset F_m\subset\ldots\subset V_n\}$. It is easy to check that the number of $G_n^0$-orbits on $X_n$ is not less than $m$. Consequently, by Theorem~\ref{theo:Omega_n_cap}, the number of $G^0$-orbits on $X$ is not less than $m$. The proof for $\SU(p,\infty)$, $p>0$, is complete.

\textsc{\underline{Case $G^0=\SU(0,\infty)$}}. Evident.

\textsc{\underline{Case $G^0=\SL(\infty,\Rp)$}}. First, let $X=\Grr{F}{V}$ for a subspace $F\subset V$ compatible with $E$. If $\dim F=k<\infty$, then $X$ consists of all $k$-dimensional subspaces of $V$. We claim that the number of $G^0$-orbits on $X$ equals $k+1$. Indeed, pick $A$, $B\in X$ and $n\geq k+1$ such that $A$, $B\in\iota_n(X_n)$ (recall that $\dim V_n=2n-1$). Clearly, if $A$ and $B$ belong to the same $G^0$-orbit, then
\begin{equation}
\dim A\cap\tau(A)=\dim B\cap\tau(B).\label{formula:signature_SL_n_R}
\end{equation}
Since $n\geq k+1$ and $V_n$ is $\tau$-stable, $\dim A\cap\tau(A)$ can be an arbitrary integer number from $0$ to $k$, hence the number of $G^0$-orbits on $X$ is at least $k+1$.

On the other hand, suppose that (\ref{formula:signature_SL_n_R}) is satisfied. Let $A'$, $B'$ be complex subspaces of $A$, $B$ respectively such that $A=A'\oplus(A\cap\tau(A))$ and $B=B'\oplus(B\cap\tau(B))$. Clearly, $A'\cap\tau(A')=B'\cap\tau(B')=0$. Furthermore, it is easy to see that
\begin{equation*}
\begin{split}
A+\tau(A)&=(A\cap\tau(A))\oplus(A'\oplus\tau(A')),\\
B+\tau(B)&=(B\cap\tau(B))\oplus(B'\oplus\tau(B')).
\end{split}
\end{equation*}
For simplicity, set $A_{\tau}=A+\tau(A)$, $A'_{\tau}=A'\oplus\tau(A')$, $A^{\tau}=A\cap\tau(A)$, and define $B_{\tau}$, $B'_{\tau}$, $B^{\tau}$ similarly. Then $A_{\tau}=A^{\tau}\oplus A'_{\tau}$, $B_{\tau}=B^{\tau}\oplus B'_{\tau}$. Note that all these subspaces are defined over $\Rp$. By \cite[Lemma 2.1]{HuckleberrySimon1}, the $\SL(A_{\tau}',\Rp)$-orbit of $A'$ is open in the corresponding grassmannian. Further\-more, there are two open $\SL(A_{\tau}',\Rp)$-orbits on this grassmannian, and their union is a single $\GL(A_{\tau}',\Rp)$-orbit. Hence there exists an operator $\psi\colon A_{\tau}\to B_{\tau}$ which is defined over $\Rp$ and maps $A'_{\tau}$, $A'$, $A^{\tau}$ to $B'_{\tau}$, $B'$, $B^{\tau}$ respectively. Since $A_{\tau}$ and $B_{\tau}$ are defined over $\Rp$ (i.e., are $\tau$-invariant), there exist\break $\tau$-invariant complements $A_0$, $B_0$ of $A_{\tau}$, $B_{\tau}$ in $V_n$. Thus one can extend $\psi$ to an operator $\nu\in\GL(V_n,\Rp)$ such that $\nu(A)=B$. Finally, since $\dim V_n$ is odd, we can scale $\nu$ to obtain an operator from $\SL(V_n,\Rp)$ which maps $A$ to $B$, as required.

At the contrary, assume that $\dim F=\infty$. As it was shown above, given $m\geq n$, two finite-dimensional spaces $A$, $B\in V_n$ belong to the same $G_m^0$-orbit if and only if $\dim A\cap\tau(A)=\dim B\cap\tau(B)$, so the number of $G_m^0$-orbits on the grassmannian of $k_n$-subspaces of $V_n$ equals $k_n+1$. But we have $\lim_{n\to\infty}k_n=\infty$, so the number of $G^0$-orbits on $X$ is infinite by Theorem~\ref{theo:Omega_n_cap}.

Now, consider the general case $X=\Fl(\Fo,E)$. We claim that, given a type\break $d=(d_1,\ldots,d_r)$, there exists a number $u(d)$ such that the number of $G_n^0$-orbits on the flag variety $\Fl(d, V_n)$ is less or equal than~$u(d)$, i.e., this upper bound depends only on $d$, but not on the dimension of $V_n$. To prove this, denote by $K_n=\SO(V_n)$ the subgroup of $G_n$ preserving the bilinear form $$\beta_n(x,y)=\sum_{i=1}^{2n-1}x_iy_i,~x=\sum_{i=1}^{2n-1}x_ie_i,~y=\sum_{i=1}^{2n-1}y_i\in V_n.$$ By Matsuki duality \cite{BremiganLorch1}, there exists a one-to-one correspondence between the set of $K_n$-orbits and the set of $G_n^0$-orbits on $\Fl(d,V_n)$. Hence our claim follows immediately from (\ref{formula:Witt_for_flags}), because \cite[Theorem~6.2]{Huang1} holds for nondegenerate symmetric bilinear forms.

Finally, suppose that $\Fo$ is of finite type. Let $\Au$, $\Bu$, $N$ be as for $\SU(\infty,\infty)$. Note that the form $\beta_n$ is nondegenerate, hence the $\beta_n$-orthogonal complement to a subspace $B\subset V_n$ is of dimension $\codim_{V_n}B$. Arguing as for $\SU(\infty,\infty)$ and applying our remark about Matsuki duality, we conclude that there exists a number $u(\Fo)$ such that the number of $G_n^0$-orbits on $X_n$ is less or equal to $u(\Fo)$ for every $n\geq N$. It follows from Theorem ~\ref{theo:Omega_n_cap} that the total number of $G^0$-orbits on $X$ is also less or equal to $u(\Fo)$. Finally, if~$\Fo$ is not of finite type, then, as in the case of $\SU(\infty,\infty)$, one can use $G$-equivariant projections from $X$ onto ind-grassmannians to show that the number of $G^0$-orbits on the ind-variety $X$ is infinite. The proof for $G^0=\SL(\infty,\Rp)$ is complete.

\textsc{\underline{Case $G^0=\SL(\infty,\Hp)$}}. Denote by $\kappa_n$ an antisymmetric bilinear form on $V_n$ defined by $$\kappa_n(e_{2i-1},e_{2i})=1,~\kappa_n(e_{2i},e_{2i-1})=-1,~\kappa_n(e_i,e_j)=0\text{ for }|i-j|>1.$$ Let $K_n$ be the subgroup of $G_n$ preserving this form. Then $K_n\cap G_n^0$ is a maximal compact subgroup of $G_n^0$ (see, e.g., \cite{FelsHuckleberryWolf1}), so, by duality, given $d$, there exists a bijection between the set of $K_n$-orbits and the set of $G_n^0$-orbits on the flag variety $\Fl(d,V_n)$. Since $K_n$ is isomorphic to $\Sp_{2n}(\Cp)$, we can argue as for $\SL(\infty,\Rp)$ to complete the proof.}

\exam{Let $X=\Grr{k}{V}$ for $k<\infty$. Then
\begin{equation*}
|X/G^0|=\begin{cases}
(p+1)(p+2)/2&\text{for }G^0=\SU(p,\infty),~p\leq k,\\
(k+1)(k+2)/2&\text{for }G^0=\SU(p,\infty),~k\leq p,\\
(k+1)(k+2)/2&\text{for }G^0=\SU(\infty,\infty),\\
k+1&\text{for }G^0=\SL(\infty,\Rp),\\
[k/2]+1&\text{for }G^0=\SL(\infty,\Hp).\\
\end{cases}
\end{equation*}
For $\SU(p,\infty)$ and $\SU(\infty,\infty)$, this follows from the formula for the number of $\SU(p,n-p)$-orbits on a finite-dimensional grassmannian, see Subsection~\ref{sst:finite_dim_case}. For $\SL(\infty,\Rp)$, this was proved in Proposition~\ref{prop:finiteness}; the proof for $\SL(\infty,\Hp)$ is similar to the case of $\SL(\infty,\Rp)$.}

As a corollary of Theorem~\ref{theo:Omega_n_cap}, we describe the degeneracy order on the set $X/G^0$ of $G^0$-orbits on an arbitrary ind-variety $X=\Fl(F,E)$ of finite type. By definition, $\Omega\leq\Omega'\iff\Omega\subseteq\overline{\Omega}'$. We define the partial order on the set $X_n/G_n^0$ of $G_n^0$-orbits on $X_n$ in a similar way.

\corop{Suppose\label{coro:degen_order} the number of $G^0$-orbits on $X=\Fl(\Fo,E)$ is finite. Then there exists $N$ such that $X/G^0$ is isomorphic as partially ordered set to $X_n/G_n^0$ for each $n\geq N$.}{Given a $G^0$-orbit $\Omega$ on $X$, there exists $n$ such that $\Omega\cap\iota_n(X_n)$ is nonempty. Since there are finitely many $G^0$-orbits on $X$, there exists $N$ such that $\Omega\cap\iota_N(X_N)$ is nonempty for each orbit~$\Omega$. By Theorem~\ref{theo:Omega_n_cap}, given $n\geq N$ and a $G^0$-orbit $\Omega$ on $X$, there exists a unique $G_n^0$-orbit $\Omega_n$ on $X_n$ such that $\iota_n^{-1}(\Omega\cap\iota_n(X_n))=\Omega_n$. Hence, the map $$\alpha_n\colon X/G^0\to X_n/G_n^0,~\Omega\mapsto \Omega_n$$ is well defined for each $n\geq N$. It is clear that this map is bijective. It remains to note that, by the definition of the topology on $X$, a $G^0$-orbit $\Omega$ is contained in the closure of a $G^0$-orbit $\Omega'$ if and only if $\Omega_n$ is contained in the closure of $\Omega_n'$ for all $n\geq N$. Thus, $\alpha_n$ is in fact an isomorphism of the partially ordered sets $X/G^0$ and $X_n/G_n^0$ for each $n\geq N$.}

\section{Open and closed orbits}\fakesect

In this section we provide necessary and sufficient conditions for a given $G^0$-orbit on $X=\Fl(\Fo,E)$ to be open or closed. We also prove that $X$ has both an open and a closed orbit if and only if the number of orbits is finite.

First, consider the case of open orbits. Pick any $n$. Recall \cite{HuckleberryWolf1}, \cite{Wolf2} that the $G_n^0$-obit of a flag $\Au=\{A_1\subset A_k\subset\ldots\subset A_k\}\in X_n$ is open if and only if
\begin{equation*}
\begin{split}
&\text{for\label{formula:nondegen_fin_dim} $G^0=\SU(p,\infty)$ or $\SU(\infty,\infty)$: all $A_i$'s are nondegenerate with respect to $\omega$;}\\
&\text{for $G^0=\SL(\infty,\Rp)$: for all $i$, $j$, $\dim A_i\cap\tau(A_j)$ is minimal},\\
&\hphantom{\text{for $G^0=\SL(\infty,\Rp)$: }}\text{i.e., equals $\max\{\dim A_i+\dim A_j-\dim V_n,~0\}$;}\\
&\text{for $G^0=\SL(\infty,\Hp)$: for all $i$, $j$, $\dim A_i\cap J(A_j)$ is minimal in the above sense}.\\
\end{split}
\end{equation*}

Note that, for any two generalized flags $\Fo_1$ and $\Fo_2$ in $X$, there is a canonical identification of $\Fo_1$ and $\Fo_2$ as linearly ordered sets. For a space $F\in\Fo_1$, we call the image of $F$ under this identification the space in $\Fo_2$ \emph{corresponding} to $F$.

Fix an antilinear operator $\mu$ on $V$. A point $\Go$ in $X=\Fl(\Fo,E)$ is \emph{in general position with respect to} $\mu$ if $F\cap\mu(H)$ does not properly contain $\wt F\cap\mu(\wt H)$ for all $F,~H\in\Go$ and all $\wt\Go\in X$, where $\wt F$, $\wt H$ are the spaces in $\wt\Go$ corresponding to $F$, $H$ respectively. A similar definition can be given for flags in~$X_n$. Note that, for $G^0=\SL(\infty,\Rp)$ or $\SL(\infty,\Hp)$, the $G_n^0$-orbit of $\Au\in X_n$ is open if and only if $\Au$ is in general position with respect to $\tau$ or~$J$ respectively.

With the finite-dimensional case in mind, we give the following
\defi{A generalized\label{defi:gen_position} flag $\Go$ is \emph{nondegenerate} if
\begin{equation*}
\begin{split}
\text{for }G^0&=\SU(p,\infty)\text{ or $\SU(\infty,\infty)$:}\\
&\text{each $F\in\Go$ is nondegenerate with respect to $\omega$;}\\
\text{for }G^0&=\SL(\infty,\Rp)\text{ or $\SL(\infty,\Hp)$}:\\
&\text{$\Go$ is in general position with respect to $\tau$ or $J$ respectively.}
\end{split}
\end{equation*}}
\vspace{-0.5cm}

\nota{A generalized flag being nondegenerate with respect to $\omega$ can be thought of as being ``in general position with respect to $\omega$''. Therefore, all conditions in Definition~\ref{defi:gen_position} are clearly analogous.}
\propp{The $G^0$-orbit $\Omega$ of $\Go\in X$ is open if and only if $\Go$ is nondegenerate.}{By the definition of the topology on $X$, $\Omega$ is open if and only if\break $\Omega_n=\iota_n^{-1}(\Omega\cap\iota_n(X_n))$ is open for each $n$.

First, suppose $G^0=\SU(p,\infty)$ or $\SU(\infty,\infty)$. To prove the claim in this case, it suffices to show that $A\in\Go$ is nondegenerate with respect to $\omega$ if and only if $\restr{\omega}{A\cap V_n}$ is nondegenerate for all $n$ for which $m_n\geq m_{\Go}$. This is straightforward. Indeed, if $A$ is degenerate, then there exists $v\in A$ such that $\omega(v,w)=0$ for all $w\in A$. Let $v\in V_{n_0}$ for some $n_0$ with $m_{n_0}\geq m_{\Go}$. Then $\restr{\omega}{A\cap V_{n_0}}$ is degenerate. On the other hand, if $v\in A\cap V_n$ is orthogonal to all $w\in A\cap V_n$ for some $n$ such that $m_n\geq m_{\Go}$, then $v$ is orthogonal to all $w\in A$ because $e_m$ is orthogonal to $V_n$ for $e_m\notin V_n$. The result follows.

Second, consider the case $G^0=\SL(\infty,\Rp)$. Suppose $\Omega$ is open, so $\Omega_n$ is open for each $n$ satisfying $m_n\geq m_{\Go}$. Assume $\Go\in X$ is not nondegenerate. Then there exist $\wt\Go\in X$ and $A,~B\in\Go$ such that $\wt A\cap\tau(\wt B)\subsetneq A\cap\tau(B)$, where $\wt A,~\wt B$ are the subspaces in $\wt\Go$ corresponding to $A,~B$ respectively. Let $v\in (A\cap\tau(B))\setminus(\wt A\cap\tau(\wt B))$, and $n$ be such that $v\in V_n$. Since $V_n$ is $\tau$-invariant, we have $v\in(A_n\cap\tau(B_n))\setminus(\wt A_n\cap\tau(\wt B_n))$ where $A_n=A\cap V_n$, $B_n=B\cap V_n$, $\wt A_n=\wt A\cap V_n$, $\wt B_n=\wt B\cap V_n$. This means that $\wt A_n\cap\tau(\wt B_n)$ is properly contained in $A_n\cap\tau(B_n)$. Hence, $\Go_n$ is not in general position with respect to $\restr{\tau}{V_n}$, which contradicts the condition that $\Omega_n$ is open.

Now, assume that $\Omega_n$ is not open for some $n$ with $m_n\geq m_{\Go}$. This means that there exist $A_n,~B_n\in\Go_n=\iota_n^{-1}(\Go)$ and $\wt\Go_n\in X_n$ so that $A_n\cap\tau(B_n)$ properly contains $\wt A_n\cap\tau(\wt B_n)$, where $\wt A_n$ and $\wt B_n$ are the respective subspaces in~$\wt\Go_n$ corresponding to $A_n$ and $B_n$. Since $\tau(e_{n+1})=e_{n+1}$, the space $A_{n+1}\cap\tau(B_{n+1})$ properly contains\break $\wt A_{n+1}\cap\tau(\wt B_{n+1})$, where $A_{n+1}$, $B_{n+1}$, $\wt A_{n+1}$, $\wt B_{n+1}$ are the respective images of $A_n$, $B_n$, $\wt A_n$, $\wt B_n$ under the embedding $X_n\hookrightarrow X_{n+1}$. Repeating this procedure, we see that $\Go$ is not nondegenerate. The result follows.

The case $G^0=\SL(\infty,\Hp)$ can be considered similarly.}

We say that two generalized flags \emph{have the same type} if there is an automorphism of $V$ transforming one into the other. Clearly, two $E$-commensurable generalized flags always have the same type. On the other hand, it is clearly not true that two generalized flags having the same type are $\wt E$-commensurable for some basis $\wt E$.

It turns out that, for $G^0=\SU(p,\infty)$ and $\SU(\infty,\infty)$, the requirement for the existence of an open orbit on an ind-variety of the form $\Fl(\Fo,E)$ imposes no restriction on the type of the flag $\Fo$. More precisely, we have

\corop{If $G^0=\SU(p,\infty)$, $0\leq p<\infty$, then $X$ always\label{coro:SU_open} has an $G^0$-open orbit. If $G^0=\SU(\infty,\infty)$, then there exist a basis $\wt E$ of $V$ and a generalized flag $\wt\Fo$ such that $\Fo$ and $\wt\Fo$ are of the same type and $\wt X=\Fl(\wt\Fo,\wt E)$ has an open $G^0$-orbit.}{For $\SU(p,\infty)$, let $n\geq m_{\Fo}$ be a positive integer such that the positive index of $\restr{\omega}{V_n}$ equals $p$. Let $\Go_n\in X_b$ be a flag in $V_n$ consisting of nondegenerate subspaces (i.e., the $G_n^0$-orbit of $\Go_n$ is open in $X_n$). Denote by $g$ a linear operator from $G_n$ such that $g(\Fo_n)=\Go_n$, where $\Fo_n=\iota_n^{-1}(\Fo)\in X_n$. Then, clearly, $g(\Fo)$ belongs to $X$ and is nondegenerate. Therefore the $G^0$-orbit of $g(\Fo)$ on $X$ is open.

Now consider the case $G^0=\SU(\infty,\infty)$. Let $\wt E$ be an $\omega$-orthogonal basis of $V$. Fix a bijection $E\to\wt E$. This bijection defines an automorphism $V\to V$. Denote by $\wt\Fo$ the generalized flag consisting of the images of subspaces from $\Fo$ under this isomorphism. Then $\wt\Fo$ and $\Fo$ are of the same type, and each space in $\wt\Fo$ is nondegenerate as it is spanned by a subset of $\wt E$. Thus the $G^0$-orbit of $\wt\Fo$ on $\wt X$ is open.}

\nota{Of course, in general an ind-variety $\wt X=\Fl(\wt\Fo,\wt E)$ having an open $\SU(\infty,\infty)$-orbit does not equal a given $X=\Fl(\Fo,E)$.}

The situation is different for $G^0=\SL(\infty,\Rp)$. While an ind-grassmannian $\Grr{F}{E}$ has an open orbit if and only if either $\dim F<\infty$ or $\codim_VF<\infty$, an ind-variety of the form $\wt X=\Fl(\wt\Fo,\wt E)$, where $\wt\Fo$ has the same type as the flag $\Fo$ from Example~\ref{exam:gen_flags} (ii), cannot have an open orbit as long as the basis $\wt E$ satisfies $\tau(\wt e)=\wt e$ for all $\wt e\in\wt E$. Indeed, suppose $\wh\Fo=\{\{0\}=\wh F_0\subset\wh F_1\subset\ldots\}\in\wt X$. As we pointed out in Example~\ref{exam:gen_flags}~(ii), there exists $N$ such that $\wh F_n=\wt F_n=\langle\wt e_1,~\ldots,~\wt e_n\rangle_{\Cp}$ for $n\geq N$. Pick~$n$ so that $m_n\geq\max\{2N,~\wt m_{\wh\Fo}\}$, where $\wt m_{\wh\Fo}$ is an integer such that $\wt\Fo$ and $\wh\Fo$ are compatible with respect to bases containing $\{\wt e_i,~i\geq\wt m_{\wh\Fo}\}$. Then the flag $\wh\Fo_n=\iota_n^{-1}(\wh\Fo)\in\wt X_n$ contains the subspace $\wh F_N$ which is defined over~$\Rp$. Thus, $\wh\Fo_n$ is not in general position with respect to $\restr{\tau}{V_n}$, so the $G_n^0$-orbit of $\wh\Fo_n$ in $\wt X_n$ is not open. Consequently, the $G^0$-orbit of $\wh\Fo$ in $\wt X$ is not open.

Let now $\wt X=\Fl(\wt\Fo,\wt E)$ where $\wt\Fo$ is a generalized flag having the same type as the generalized flag $\Fo$ from Example~\ref{exam:gen_flags} (iii). Recall that $$\Fo=\{\{0\}=F_0\subset F_1\subset F_2\subset\ldots\subset F_{-2}\subset F_{-1}\subset V\},$$ where $F_i=\langle e_1,~e_3,~\ldots, e_{2i-1}\rangle_{\Cp}$, $F_{-i}=\langle\{e_j,~j\text{ odd}\}\cup\{e_{2j},~j>i\}\rangle_{\Cp}$ for $i\geq1$. We claim that $\wt X$ also cannot have an open orbit. Indeed, assume $\wh\Fo$ is $\wt E$-commensurable to $\wt\Fo$. Then $\wh\Fo$ is compatible with a basis $\wh E$ of~$V$ such that $\wh E\setminus\wt E$ is finite. This means that there exists $\wt e\in\wt E$ and a finite-dimensional subspace $F\in\wh \Fo$ with $\wt e\in F$. Now, pick $n$ so that $F\subset V_n$ and $m_n\geq\max\{2\dim F,~\wt m_{\wh\Fo}\}$. Then $F\cap\tau(F)\neq0$, so\break $\wh\Fo_n=\iota_n^{-1}(\wh\Fo)\in X_n$ is not in general position with respect to $\restr{\tau}{V_n}$.

Finally, let $G^0=\SL(\infty,\Hp)$. In this case, clearly, an~ind-grass\-man\-nian $\Grr{F}{E}$ may or may not have an open orbit. A similar argument as for $\SL(\infty,\Rp)$ shows that if $\Fo$ is as in Example~\ref{exam:gen_flags} (ii), then $\wt X$ cannot have an open orbit. Surprisingly, for $G^0=\SL(\infty,\Hp)$ and $X$ as in Example~\ref{exam:gen_flags}~(iii), $\wt X$ may have an open orbit. Consider first the case of $X=\Fl(\Fo,E)$ itself. It it easy to check that if $\dim V_n=n$ then $\Fo_n$ is in general position with respect to $\restr{J}{V_n}$ for each $n$, so the orbit of $\Fo$ is open. On the other hand, if $\wt\Fo$ and $\Fo$ have the same type and each $2n$-dimensional subspace in $\wt\Fo$ is spanned by the vectors $\wt e_1,~\wt e_2,~\wt e_5,~\wt e_6,~\ldots,~\wt e_{4n-3},~\wt e_{4n-2}$, then $\wt X$ does not have an open orbit because each generalized flag $\wt E$-commensurable to $\wt F$ contains a finite-dimensional subspace $F$ such that $F\cap J(F)\neq\{0\}$.

We now turn our attention to closed orbits. The conditions for an orbit to be closed are based on the same idea for each of the real forms, but (as was the case for open orbits) the details differ.

Suppose $G^0=\SU(\infty,\infty)$ or $\SU(p,\infty)$. We call a generalized
flag $\Go$ in $X$ \emph{pseudo-isotropic} if $F\cap H^{\perp,V}$ is not
properly contained in $\wt F\cap\wt H^{\perp,V}$ for all $F,~H\in\Go$ and
all $\wt\Go\in X$, where $\wt F,~\wt H$ are the subspaces in $\wt\Go$
corresponding to $F,~H$ respectively. A similar definition can be given
for flags in $X_n$. An isotropic generalized flag, as defined in \cite{DimitrovPenkov1}, is always pseudo-isotropic,
but the converse does not hold. In the particular case when the generalized flag $\Go$ is of the form
$\{\{0\}\subset F\subset V\}$, $\Go$ is pseudo-isotropic if and only if the kernel
of the form $\restr{\omega}{F}$ is maximal over all $E$-commensurable flags
of the form $\{\{0\}\subset\wt F\subset V\}$.

Next, suppose $G^0=\SL(\infty,\Rp)$. A generalized flag $\Go$ in $X$ is \emph{real} if $\tau(F)=F$ for all $F\in\Go$. This condition turns out to be equivalent to the following condition: $F\cap\tau(H)$ is not properly contained in $\wt F\cap\tau(\wt H)$ for all $F,~H\in\Go$ and all $\wt\Go\in X$, where $\wt F,~\wt H$ are the subspaces in $\wt\Go$ corresponding to $F,~H$ respectively.

Finally, suppose $G^0=\SL(\infty,\Hp)$. We call a generalized flag $\Go$ in $X$ \emph{pseudo-quaternionic} if $F\cap J(H)$ is not properly contained in $\wt F\cap J(\wt H)$ for all $F,~H\in\Go$ and all $\wt\Go\in X$, where $\wt F,~\wt H$ are the subspaces in $\wt\Go$ corresponding to $F,~H$ respectively. If $\Go$ is \emph{quaternionic}, i.e., if $J(F)=F$ for each $F\in\Go$, then $\Go$ is clearly pseudo-quaternionic, but the converse does not hold. If the generalized flag $\Go$ is of the form $\{\{0\}\subset F\subset V\}$, then $\Go$ is pseudo-quaternionic if and only if $\codim_F(F\cap J(F))\leq1$.

\propp{The $G^0$-orbit\label{prop:closed_orbit} $\Omega$ of $\Go\in X$ is closed if and only if
\begin{equation*}
\begin{split}
&\text{$\Go$ is pseudo-isotropic for $G^0=\SU(\infty,\infty)$ and $\SU(p,\infty)$\textup;}\\
&\text{$\Go$ is real for $G^0=\SL(\infty,\Rp)$\textup;}\\
&\text{$\Go$ is pseudo-quaternionic for $G^0=\SL(\infty,\Hp)$.}\\
\end{split}
\end{equation*}}
{First consider the finite-dimensional case, where there is a unique closed $G_n^0$-orbit on $X_n$ (see Theorem~\ref{theo:finite_dim_case}). For all real forms the conditions of the proposition applied to finite-dimensional flags in $V_n$ are easily checked to be closed conditions on points of $X_n$. Therefore, the $G_n^0$-orbit of a flag in $V_n$ is closed if and only if this flag satisfies the conditions of the proposition at the finite level.

Let $G^0=\SU(\infty,\infty)$ or $\SU(p,\infty)$. Suppose $\Omega$ is closed, so $\Omega_n$ is closed for each $n$ satisfying $m_n\geq m_{\Go}$. Assume $\Go$ is not pseudo-isotropic. Then there exist $\wt\Go\in X$ and $A,~B\in\Go$ such that $\wt A\cap\wt B^{\perp,V}\supsetneq A\cap B^{\perp,V}$, where $\wt A,~\wt B$ are the subspaces in $\wt\Go$ corresponding to $A,~B$ respectively. Let $v\in (\wt A\cap\wt B^{\perp,V})\setminus(A\cap B^{\perp,V})$, and $n$ be such that $v\in V_n$ and $m_n\geq m_{\Go}$. Then $v\in(\wt A_n\cap\wt B_n^{\perp,V_n})\setminus(A_n\cap B_n^{\perp,V_n})$, where $A_n=A\cap V_n$, $B_n=B\cap V_n$, $\wt A_n=\wt A\cap V_n$, $\wt B_n=\wt B\cap V_n$, because $B^{\perp,V}\cap V_n=B_n^{\perp,V_n}$. This means that $A_n\cap B_n^{\perp,V_n}$ is properly contained in $\wt A_n\cap\wt B_n^{\perp,V_n}$. Hence $\Go_n$ is not pseudo-isotropic, which contradicts the condition that $\Omega_n$ is closed.

Now, assume that $\Omega_n$ is not closed for some $n$ with $m_n\geq m_{\Go}$. This means that there exist $A_n,~B_n\in\Go_n=\iota_n^{-1}(\Go)$ and $\wt\Go_n\in X_n$ such that $A_n\cap B_n^{\perp,V_n}$ is properly contained in $\wt A_n\cap\wt B_n^{\perp,V_n}$, where $\wt A_n,~\wt B_n$ are the subspaces in~$\wt\Go_n$ corresponding to $A_n$, $B_n$ respectively. Since each $e\in E_{n+1}\setminus E_n$ is orthogonal to $V_n$, $A_{n+1}\cap B_{n+1}^{\perp,V_{n+1}}$ is properly contained in $\wt A_{n+1}\cap\wt B_{n+1}^{\perp,V_{n+1}}$, where $A_{n+1}$, $B_{n+1}$, $\wt A_{n+1}$, $\wt B_{n+1}$ are the respective images of $A_n$, $B_n$, $\wt A_n$, $\wt B_n$ under the embedding $X_n\hookrightarrow X_{n+1}$. Repeating this procedure, we see that $\Go$ is not pseudo-isotropic. The result follows.

Let $G^0=\SL(\infty,\Rp)$. As above, given $n$, denote $\Go_n=\iota_n^{-1}(\Go)$. Note that, given $F\in\Go$, $\tau(F)=F$ if and only if $F_n$ is defined over $\Rp$, i.e., $\tau(F_n)=F_n$ where\break $F_n=F\cap V_n$, because $V_n$ is $\tau$-invariant. The $G_n^0$-orbit $\Omega_n$ of $\Go_n$ is closed if and only if each subspace in $\Go_n$ is defined over~$\Rp$. Hence if $\tau(F)=F$ for all $F\in\Go$, then $\Omega_n$ is closed for each $n$ (so $\Omega$ is closed), and vice versa.

The proof for $G^0=\SL(\infty,\Hp)$ is similar to the case of $\SU(\infty,\infty)$ and is based on the following facts: if $A$ is a subspace of $V$, then $J(A)\cap V_n=J(A\cap V_n)$ for all $n$; the subspace $\langle E_{n+1}\setminus E_n\rangle_{\Cp}$ is $J$-invariant for all $n$.}

\corop{If $G^0=\SL(\infty,\Rp)$\textup, then $X=\Fl(\Fo,E)$ always\label{coro:SU_closed} has a closed orbit.}{The $G^0$-orbit of the generalized flag $\Fo$ is closed because $\tau(e)=e$ for all basic vectors $e\in E$.}

If $G^0=\SU(p,\infty)$ for $0\leq p<\infty$, then, for some choice of $\Fo$, $X=\Fl(\Fo,E)$ clearly may have a closed orbit. On the other hand, assume that $\Fo$ is as in Example~\ref{exam:gen_flags} (ii), and $E$ is an orthogonal basis for $\omega$. If $\wh\Fo\in X$, then then there exists $N$ such that $\wh F_n=F_n$ (and so $\wh F_n$ is nondegenerate) for $n\geq N$. But there obviously exists $\Ho\in X$ such that $H_N$ is degenerate, so $\wh\Fo$ is not pseudo-isotropic. Thus, $X$ does not contain a closed orbit.

Let $G^0=\SU(\infty,\infty)$. Obviously, the ind-grassmannian $\Grr{F}{E}$ may have or may not have a closed orbit. If $\wt X=\Fl(\wt\Fo,\wt E)$, where $\wt\Fo$ is a generalized flag having the same type as the generalized flag $\Fo$ from Example~\ref{exam:gen_flags} (ii) and $\wt E$ satisfies all required conditions, then $\wt X$ does not have a closed orbit. Indeed, assume $\wh\Fo\in\wt X$, then $\wh F_n$ contains $V_k$ for certain $n$ and $k$. The form $\restr{\omega}{V_k}$ is nondegenerate, hence $\wh F_n$ is not isotropic. There exists an isotropic subspace $I$ of $V$ of dimension $n=\dim\wh F_n$ containing $\wh F_n\cap\wh F_n^{\perp,V}$, and it is easy to see that there exists $\wh\Fo_0\in\wt X$ such that $I$ is the subspace of $\wh\Fo_0$ corresponding to $\wh F_n$. Thus, $\wh\Fo$ is not pseudo-isotropic.

Now, suppose $\Fo$ is as in Example~\ref{exam:gen_flags} (iii). Here $\wt X$ may or may have not a closed orbit. For example, assume that $\wt E$ is an $\omega$-orthogonal basis of $V$. Then each $\wh\Fo\in\wt X$ contains a nonisotropic finite-dimensional subspace, and, arguing as in the previous paragraph, we see that $\wh\Fo$ is not pseudo-isotropic. On the other hand, suppose that $e_{2i-1}=e_{2i-1}'+e_{2i}'$ and $e_{2i}=e_{2i-1}'-e_{2i}'$ for all $i$, where $\{e_1',~e_2',~\ldots\}$ is an $\omega$-orthogonal basis with $\omega(e_{2i-1}',e_{2i-1}')=-\omega(e_{2i},e_{2i})=1$. In this case, one can easily check that $\Fo$ is pseudo-isotropic, so its $G^0$-orbit in $X$ is closed.

Finally, let $G^0=\SL(\infty,\Hp)$. Here, in all three cases (i), (ii), (iii) of Example~\ref{exam:gen_flags}, if $\wt X=\Fl(\wt\Fo,\wt E)$ for a generalized flag $\wt\Fo$ having the same type as $\Fo$, then $\wt X$ may or may not have a closed orbit. Consider, for instance, case (ii). The flag $\Fo$ itself is pseudo-quaternionic, so its $G^0$-orbit in $X$ is closed. On the other hand, if each $(4n+2)$-dimensional subspace in $\wt\Fo$ is spanned by $\{e_i,~i\leq 4n\}\cup\{e_{4n+1},~e_{4n+3}\}$, then $\wt X$ does not have a closed orbit.

Combining our results on the existence of open and closed orbits, we now obtain the following corollary.

\corop{For a given real form $G^0$ of $G=\SL(\infty,\Cp)$\textup, $G^0\neq\SU(p,\infty)$\textup,\break $0<p<\infty$\textup, an ind-variety of generalized flags $X=\Fl(\Fo,E)$ has both an open and a closed $G^0$-orbits if\textup, and only if\textup, there are only finitely many $G^0$-orbits on $X$.}{If $X$ has finitely many $G^0$-orbits, then the existence of an open orbit is obvious, and the existence of a closed orbit follows immediately from Corollary~\ref{coro:degen_order}.

Assume that $X$ has both an open and a closed $G^0$-orbit. Let\break $G^0=\SU(\infty,\infty)$. Fix a nondegenerate generalized flag ${\Ho}\in X$ (lying on an open $G^0$-orbit). Suppose that there exists a subspace ${F}\in{\Ho}$ satisfying $\dim {F}=\codim_V{F}=\infty$. Since $X$ has a closed $G^0$-orbit, there exists a pseudo-isotropic generalized flag ${\wt\Ho}\in X$. Let ${\wt F}$ be the subspace in ${\wt\Ho}$ corresponding to ${F}$. Since ${\Ho}$ and ${\wt\Ho}$ are $E$-commensurable to $\Fo$, there exists $n$ such that ${F}={A}\oplus B$ and ${\wt F}={\wt A}\oplus B$, where ${A}$, ${\wt A}$ are subspaces of $V_n$ and $B$ is the span of a certain infinite subset of $E\setminus E_n$; in particular, $B$ is a subspace of $\overline V_n=\langle E\setminus E_n\rangle_{\Cp}$.

The restriction of $\omega$ to $B$  is nondegenerate, because $V_n$ and $\overline V_n$ are orthogonal. This implies that $B^{\perp,\overline V_n}\cap B=\{0\}$. But ${\wt F}^{\perp,V}={\wt A}^{\perp,V_n}\oplus B^{\perp,\overline V_n}$, hence\break ${\wt F}\cap {\wt F}^{\perp,V}={\wt A}\cap {\wt A}^{\perp,V_n}$. Clearly, if $B\neq\overline V_n$, then $B^{\perp,\overline V_n}\neq\{0\}$. In this case, there exists $v\in \overline V_n\setminus B$ contained in ${\wt F}^{\perp,V}$, and one can easily construct a generalized flag ${\wh\Ho}\in X$ such that ${\wt F}\cap {\wt F}^{\perp,V}\subsetneq {\wh F}\cap {\wh F}^{\perp,V}$, where ${\wh F}$ is the subspace in ${\wh\Ho}$ corresponding to ${\wt F}$, a contradiction. Thus, $B=\overline V_n$, but this contradicts the condition $\codim_V{\wt F}=\infty$.

We conclude that ${\Ho}=\Au\cup\Bu$, where each subspace in $\Au$ (resp., in $\Bu$) is of finite dimension (resp., of finite codimension). Assume that $\Fo$ is not finite, then at least one of the generalized flags $\Au$ and $\Bu$ is infinite. Suppose $\Au$ is infinite. (The case when $\Bu$ is infinite can be considered using the map $U\mapsto U^{\#}$.) Let $n$  be such that ${\Ho}$ and ${\wt\Ho}$ are compatible with bases containing $E\setminus E_n$. Let ${F}$ be a subspace in $\Au$ such that ${F}$ does not belong to~$V_n$. Then, arguing as in the previous paragraph, one can show that ${\wt\Ho}$ cannot be pseudo-isotropic, a contradiction.

Now, let $G^0=\SL(\infty,\Rp)$. Suppose that ${\Ho}\in X$ is in general position with respect to $\tau$, and ${\wt\Ho}\in X$ is real. As above, pick $n$ so that ${\Ho}$ and ${\wt\Ho}$ are compatible with bases of $V$ containing $E\setminus E_n$. Suppose for a moment that there exists a subspace ${F}\in{\Ho}$ such that ${F}\not\subset V_n$, then ${F}={A}\oplus B$, where ${A}$ is a subspace of $V_n$, and $B$ is a nonzero subspace of $\overline V_n$ spanned by a subset of $E\setminus E_n$. Similarly, the corresponding subspace ${\wt F}\in{\wt\Ho}$ has the form ${\wt F}={\wt A}\oplus B$, where $\tau({\wt A})={\wt A}$ and $\tau(B)=B$. Suppose also that $B\neq\overline V_n$, then there exist $e\in E\cap B$ and $e'\in(E\setminus E_n)\setminus B$. Let $B'\subset\overline V_n$ be spanned by $((E\cap B)\setminus\{e\})\cup\{e+ie'\}$. It is easy to check that there exists ${\wh\Ho}\in X$ such that the subspace ${\wh F}\in{\wh\Ho}$ corresponding to ${F}$ has the form ${A}\oplus B'$. Thus, ${F}\cap\tau({F})$ properly contains ${\wh F}\cap\tau({\wh F})$, a contradiction. It remains to note that if $\Fo$ is not of finite type, then such a subspace ${F}$ always exists (if necessary, after applying the map $U\mapsto U^{\#}$).

Finally, let $G^0=\SL(\infty,\Hp)$. Suppose that ${\Ho}\in X$ is in general position with respect to $J$, and ${\wt\Ho}\in X$ is pseudo-quaternionic. As above, pick $n$ so that ${\Ho}$ and ${\wt\Ho}$ are compatible with bases of $V$ containing $E\setminus E_n$. Suppose for a moment that there exists a subspace ${F}\in{\Ho}$ such that ${F}\not\subset V_n$, then ${F}={A}\oplus B$, where ${A}$ is a subspace of $V_n$, and $B$ is a nonzero subspace of $\overline V_n=\langle E\setminus E_n\rangle_{\Cp}$ spanned by a subset of $E\setminus E_n$. The corresponding subspace ${\wt F}\in{\wt\Ho}$ has the form ${\wt F}={\wt A}\oplus B$, where ${\wt A}$ is a subspace of~$V_n$. Suppose also that $\dim B\geq2$ and $\codim_{\overline V_n}B\geq2$. There exist a subspace $B'\subset\overline V_n$ and ${\wh\Ho}\in X$ such that the subspace ${\wh F}\in{\wh\Ho}$ corresponding to ${F}$ has the form ${A}\oplus B'$, and $B'\cap J(B')$ is either properly contains or is properly contained in $B\cap J(B)$. Thus, either ${\Ho}$ is not in general position with respect to $J$, or ${\wt\Ho}$ is not pseudo-quaternionic, a contradiction. It remains to note that if $\Fo$ is not of finite type, then such a subspace ${F}$ always exists (possibly, after applying the map $U\mapsto U^{\#}$).}

\begin{acknowledgement}
The first author was supported in part by the Russian Foundation for Basic Research through grants no. 14--01--97017 and 16--01--00154, by the Dynasty Foundation and by the Ministry of Science and Education of the Russian Federation, project no. 204. A part of this work was done at the Oberwolfach Research Institute for Mathematics (program ``Oberwolfach Leibniz Fellows'') and at Jacobs University Bremen, and the first author thanks these institutions for their hospitality. The second and third authors thank Professor V.K. Dobrev for the invitation to speak at the XI International Workshop ``Lie Theory and its Applications in Physics'' in Varna, 15--21 June 2015. The second author acknowledges continued partial support by the DFG through Priority Program SPP 1388 and grant PE 980/6--1.  The third author acknowledges partial support from the Dickson Emeriti Professorship at the University of California and from a Simons Foundation Collaboration Grant for Mathematicians.
\end{acknowledgement}

\end{document}